\documentclass[a4paper,twoside]{article}
\usepackage{amssymb}
\usepackage[latin1]{inputenc}
\usepackage[T1]{fontenc}
\usepackage[final]{lpretex}
\usepackage[all]{xy}
\usepackage{amscd}
\usepackage{title}

\EndOfDump

\usepackage{htmlref}

\DocVersion[Sun May 18 2003]{1}

\newcommand\equaliser[2]{{\mapright{#1}\atop\mapright{#2}}}
\DeclareRobustCommand\Sch{\mathbf{S}}
\renewcommand\Hom{\symb{Hom}}
\newcommand\Mj{\Cal{S}^\mu}
\newcommand\Si{\Cal{S}^{\iso}}

\newcommand\nq[1]{\mathbf{N}\Cal{Q}{#1}}
\newcommand\bh[1]{\Cal{B}{#1}}
\newcommand\qh[1]{\Cal{Q}{#1}}
\renewcommand\dh[1]{\Cal{D}{#1}}
\newcommand\dMj[1]{\Cal{D}^{#1}{\Mj}}

\newcommand\Mod{\symb{Mod}}

\begin{document}

\title{Strict polynomial functors and multisets}
\author{Torsten Ekedahl}
\address{Department of Mathematics\\
Stockholm University\\
SE-106 91 Stockholm\\
SWEDEN}
\email{teke@math.su.se}
\author{Pelle Salomonsson}
\address{Department of Mathematics\\
Stockholm University\\
SE-106 91 Stockholm\\
SWEDEN}
\email{p.s@tiscali.se}

\begin{abstract}
We prove a generalisation to any characteristic of a result of Macdonald that
describes strict polynomial functors in characteristic zero in terms of
representations of the groupoid of finite sets and bijections. Our result will
give an analogous description in terms of finite multisets and an extension of the
notion of bijection appropriate for multisets. A projected application is to the
description of strict polynomial monads that will give a notion generalising
(linear algebra) operads.
\end{abstract}
\maketitle

Over fields (or more generally (commutative) rings) of characteristic $0$ a
result of Macdonald (cf.\ \cite[App.\ I:A:5.3]{macdonald95::symmet+hall}) gives
an equivalence between the category of strict polynomial functors (i.e.,
functors $F$ from finitely generated projective modules to arbitrary modules for
which the structure maps $\Hom(P,Q) \to \Hom(FP,FQ)$ have been given the
structure of polynomial maps) to the category of $\Sigma$-modules; collections
$\{M_n\}$ of $\Sigma_n$-modules for each $n \ge 0$. This relation with the
symmetric groups provides the study of the category of strict polynomial
functors with a heavily combinatorial flavour which becomes even more pronounced
if one instead of $\Sigma$-modules considers the equivalent category of
representations of the groupoid of finite sets and bijections. It is the purpose
of the current article to give an analogous description of the category of
strict polynomial functors over any commutative ring. Very roughly speaking this
description is obtained by replacing finite sets with finite \emph{multisets}
and bijections with a notion that we shall call \emph{multijections}. (Note that
this is not the only ``algebraic description'' of strict polynomial functors;
they also correspond to collections of modules over Schur algebras. Our
contribution is to give a description close in spirit to Macdonald's.)

Even though there is a strong analogy between sets/bijections and
multisets/multijections there is one important conceptual difference in that
bijections are always invertible whereas multijections are not (they do have the
property that a multijection endomorphism is an isomorphism but there are
multijections also between some non-isomorphic multisets). This has as a
consequence that in the category of representations of multisets that will turn
out to be equivalent to the category of strict polynomial functors one will have
for a multijection $f$ two maps $f_*$ and $f^*$ acting on the representation in
question. These will be inverses to each other when $f$ is invertible but will
not be in general. There will be certain relations between their composites most
of which serve to make the representations be special cases of the notion of
Mackey functors introduced by Dress (cf.\ \cite{dress73::contr}). There is a
small subtlety in that the category of finite multisets and muiltijections is
not a based category in Dress' sense but only gives rise to one; we shall
describe the precise relationship. Furthermore, these functor will fulfill one
further property: Each multijection has associated to it a degree, a positive
integer and we will also have the relation $f_*f^* = \deg f$; we shall call such
functors \emph{Hecke-Mackey functors} and our result is that the category of
strict polynomial functors is equivalent to the category of Hecke-Mackey
functors on the category of multisets and multijections.

We give the precise relation with the Macdonald description by showing that over
a $\Q$-algebra our notion is equivalent to that of $\Sigma$-module, a collection
of $\Sigma_n$-modules. From our point of view it is more natural to say that one
may restrict attention to multisets all of whose elements have multiplicity one,
i.e., ordinary sets. The same idea of proof also shows that in a $p$-local
situation one may restrict attention to sets all of whose multiplicities are
powers of $p$.

We feel that that our description is interesting from a purely
representation-theoretic view point and that it would for instance be
interesting to compare our description with the approach using Schur algebras
(which has been much further developed). Partitions, which to us will simply be
a special kind of multisets, appear very explicitly in the theory of Schur
algebras and it is not clear if these appearances can be related to the r\^ole
they play in our theory. We would, however, like to point out that the main
application we have in mind goes in a different direction. In a sequel
we shall use the description of this paper to introduce the notion of Schur
operad which will then be proven to be equivalent to that of strict polynomial
monads. This gives nothing new in characteristic $0$ but will allow us to give
an operadic description of for instance $p$-Lie algebras which is not covered by
ordinary operads. These ideas will also be seen to give (yet another) algebraic
framework for describing the homology operations on infinite loop spaces.
\begin{conventions}
Throughout this article we shall let $R$ be a commutative Noetherian ring.
\end{conventions}

\begin{section}{Preliminaries}

\begin{subsection}{Multisets}

To us a \Definition{multiset} shall be a function from the class of all sets to
$\N$, the set of non-negative numbers, which is zero outside of a subset. The
value of the function on a set will be called the \Definition{multiplicity} of
the multiset at the set. We shall say that the multiset is \Definition{supported
in} a set $S$ if it is zero outside it. The \Definition{support} of a multiset
is the set where it is non-zero, we shall then also say that the multiset is
\Definition{supported on} its support. If $S$ is a multiset we shall denote its
support by $\underline{S}$ and then we shall denote the restriction of $S$ to
$\underline{S}$ by $\mu_S$. A multiset is completely determined by $\mu_S$ and
conversely any $\N$-valued function on a set can be extended by zero outside the
set to give rise to a multiset. By abuse of language we shall consistently
confuse these two descriptions (it is very convenient to have the function
defined on all sets which is why we need two descriptions). We shall also
consider each set $S$ as a multiset by the condition that the support of $\mu_S$
is $S$ and that $\mu_S$ is identically $1$ on $S$. In particular, using the
$\mu_S$ point of view a multiset is also a set (as $\mu_S$ is) and so we will
have no problem with considering multisets of multisets. Furthermore, a multiset
$T$ is a \Definition{subset} of a multiset $S$ if $T(s) \le S(s)$ for all sets
$s$. We define the \Definition{product}, $S\times T$, by the condition that
$(S\times T)(w)$ is zero unless $w$ has the form $w=(s,t)$ in which case
$(S\times T)(w):=S(s)T(t)$. Similarly, we define the \Definition{disjoint union}
of two multisets $S$ and $T$ as the multiset whose support is the disjoint union
of the supports and whose multiplicity on the copy of $\underline{S}$ is given
by the multiplicity function of $S$ and similarly for the copy of
$\underline{T}$.

We shall mainly be interested in \Definition{finite multisets}, that is, multisets
with finite support and for such a multiset $S$ its \Definition{cardinality} is
defined to be $|S| := \sum_{s\in\underline{S}}\mu_S(s)$. We shall also use
standard set-theoretic notation for multisets in cases where the meaning can be
easily inferred. So we may for instance express the cardinality of a multiset
$S$ as $\sum_{s\in S}1$ where it is hence understood that the sum is summed with
multiplicities. More generally for a finite multiset $S$ of elements of a
commutative monoid we mean by $\sum_{s \in S}s$ the sum $\sum_{t \in
\underline{S}}\mu_S(t)t$.

If $S$ and $T$ are multisets we define their \Definition{sum} by
$(S+T)(x):=S(x)+T(x)$. This makes the set of multisets supported in a fixed set
an abelian monoid (more precisely it is the power of $\N$ over the set).

If $S$ and $T$ are finite multisets then a \Definition{multijection} from $S$ to
$T$ is a map \map{f}{\underline{S}}{\underline{T}} such that
$\mu_T(t)=\sum_{f(s)=t}\mu_S(s)$ for $t \in \underline{T}$. It is clear that if
\map{g}{\underline{T}}{\underline{U}} is the underlying map of another
multijection then the composite $g\circ f$ is the underlying map of a
multijection which we shall call the \Definition{composite} of $f$ and $g$. From
this it is equally clear that the finite multisets and the multijections between
them form a category, $\Mj$. Note that the multijection as well as the multiset
$T$ is determined by the map $f$ and conversely, given a finite multiset $S$ and
a map of sets \map{f}{\underline{S}}{T} we can define $\mu_T$ by $\mu_T(t)=
\sum_{s \in f^{-1}(t)}\mu_S(s)$ and then $f$ becomes a multijection. We shall
call this multijection the \Definition{induced multijection}. Note further that
$\Mj$ has a symmetric monoidal structure; the disjoint union of multijections is
again a multijection. A multijection $S \to T$ will be called \Definition{final}
if the support of $T$ consists of a single point. It is then clear that every
multijection decomposes, up to isomorphism and in an essentially unique way, as
a disjoint union of final multijections.

Isomorphism classes of finite multisets correspond exactly to unordered
partitions. We shall now describe an association of a specific representative to
each unordered partition. If $[1^{n_1},\dots,k^{n_k}]$ is an unordered partition
then we associate to it the multiset with support in the positive integers that
gives $1,2,\dots,n_1$ multiplicity $1$, $n_1+1,\dots,n_1+n_2$ multiplicity $2$
etc. We shall adopt the custom of similar situations and in fact let
$[1^{n_1},\dots,k^{n_k}]$ mean both the partition and this multiset
representative. (Note that thus $[n]$ means $\{1,1,1,\dots,1\}$ which differs
from a not uncommon notation where it is used for $\{1,2,\dots,n\}$. We will
instead have $[1^n]=\{1,2,\dots,n\}$.)

We shall also need to consider ordered partitions that we shall identify
with finite multisets supported in the non-negative integers. The notation
$[n_1,n_2,\dots,n_k]$, $n_i$ non-negative integers, will sometimes be used for
the partition which has multiplicity $n_i$ at $i$, $1\le i \le k$ and zero
otherwise. (This gives a slight notational clash in that for instance $[n]$
could denote both an unordered and an ordered partition, this should not however
cause any confusion.)
\begin{proposition}\label{generators}
Every multijection can be written as a composite of isomorphisms of multisets
and the disjoint union of a final morphism of the type $[n_1,n_2] \to
[n_1+n_2]$ with identity maps.
\begin{proof}
This is clear as every multiset is isomorphic to an ordered partition and every
multijection can be written as the composite of maps $\alpha \to \beta$ that are
bijections outside of one point of $\underline{\beta}$ and those in the latter
category is a composition of such, such that the cardinality of the fibres are $1$ or
$2$.
\end{proof}
\end{proposition}
A \Definition{multi-map} from a finite multiset $\alpha$ to a (multi)set $S$ is a
submultiset $\Gamma$ of $\alpha\times S$ such that projection on the first
factor is a multijection. Just as for functions we shall sometimes make a
distinction between the multi-map and its \Definition[graph of multi-map]{graph}
which is $f$ considered as a submultiset of $\alpha\times S$.
\begin{remark}
This definition is in analogy with the definition of a function on sets as a
subset of the product such that projection on the first factor is a
bijection. Intuitively, if a multiset can be thought of as a set with elements
repeated according to their multiplicity, then a multi-map may take different
values on different copies of the same element. The multijection which is the
projection on the first factor is then the smallest multijection needed to make
this function one-valued on the underlying set. A map on the other hand does not
allow this. Note furthermore that there is no natural way of composing
multi-maps.
\end{remark}
A multi-map is a \Definition{map} if the projection on the first factor is
actually an isomorphism. A multi-map is a \Definition{bimultijection} if the
projection on the second factor is also a multijection. Note that multi-maps
$\alpha \to S$ may be identified with isomorphism classes of diagrams $\alpha \leftarrow
\gamma \to S$, where $\alpha \leftarrow \gamma$ is a multijection and $\gamma
\to S$ is a map whose underlying function $\underline{\gamma} \to \underline{S}$
is injective on the fibres of $\underline{\alpha} \leftarrow \underline{\gamma}$.

We shall need the following generalisation of shuffles to the multiset situation.
\begin{definition}
If \map{f}{\beta}{\Mj} is a map from a finite multiset to the class of finite
multisets and $\alpha := \sum_{s \in \beta}f(s)$ then the
\Definition{multi-shuffle} associated to $f$ is the induced multijection
$\disjunion_{s \in \beta}f(s) \to \alpha$.\footnote{Note that according to the
general principles $f(s)$ will appear $\mu_\beta(s)$ times in the disjoint
union.} If we want to be more specific we shall call it a
\Definition{$\beta$-heaped multi-shuffle}.
\end{definition}
\end{subsection}
\begin{subsection}{Hecke-Mackey functors}

We recall (cf., \cite[\S4]{dress73::contr}) the definition of a
\Definition{based category} making the trivial modification that we only require
that there be a finite number of isomorphism classes in each component.
\begin{itemize}
\item[(M1)] $\cA$ has finite limits and finite coproducts,

\item[(M2)] the two squares of a diagram
\begin{displaymath}
\begin{array}{ccccc}
X' & \longrightarrow & Z' & \longleftarrow & Y'\\
\downarrow &&      \downarrow && \downarrow\\
X& \longrightarrow & X\coprod Y & \longleftarrow & Y
\end{array}
\end{displaymath}
are cartesian precisely when the top row is a coproduct,

\item[(M3)] there is a finite set of isomorphism classes of indecomposable
objects in each component of $\cA$ and every object of $\cA$ is a coproduct of
indecomposable objects, and

\item[(M4)] for indecomposable $Z$ and $Z'$, $\Hom_{\cA}(Z,Z')$ is finite and an
endomorphism of an indecomposable object is an isomorphism.
\end{itemize}
As is essentially noticed in \cite{dress73::contr} a based category is
determined by its full subcategory of indecomposable objects. We shall need to
make this explicit as well as giving a characterisation of that subcategory. For
that recall (cf., \cite[Déf.~5.2]{baues01::fonct+mackey}) that a category has
\Definition{weak fibre products} if for every pair of morphisms with the same
target there are a finite number of commutative squares with the left and bottom
maps equal to the two given ones and such that any other square with this
property factors uniquely through precisely one of them.
\begin{definition}
An \Definition{atomic category} is a category $\cA$ with the following properties.
\begin{itemize}
\item[(A1)] $\cA$ has weak fibre products and a final object and no initial object,

\item[(A3)] there is a finite number of isomorphism classes of objects in each
component of $\cA$, and

\item[(A4)] $\Hom_{\cA}(Z,Z')$ is finite for any objects $Z$ and $Z'$ of $\cA$
and an endomorphism of an object is an isomorphism.
\end{itemize}
\end{definition}
We shall now see that these two notions are essentially the same.
\begin{proposition}
The full subcategory of indecomposable objects of a based category form an
atomic category and for each atomic category there is a, unique up to
equivalence, based category whose full subcategory of indecomposable objects is
equivalent to the given atomic category.
\begin{proof}
Start with a based category.  By (M1) (cf., last paragraph of
\cite[p.~202]{dress73::contr}) any map from an indecomposable object to a
coproduct factors through one of the factors. From this the axioms for an atomic
category are immediate. Furthermore, again from (M1), if $X_1,\dots,X_n$ and $Y_1,\dots,Y_m$
are indecomposable then 
\begin{displaymath}
\Hom(\coprod_iX_i,\coprod_jX_j) = \coprod_j\prod_i\Hom(X_i,Y_j)
\end{displaymath}
which shows that the category is determined up to equivalence by the full
subcategory of indecomposable objects as well as giving a recipe for
constructing a based category given an atomic; its objects are sequences
(possibly empty to get an initial object) $(X_1,\dots,X_n)$ and morphisms are
given by the above formula (alternatively one gets an equivalent category by
using finite multisets of objects) and one checks easily that one gets in this
way a based category which contains a full subcategory equivalent to the atomic
one consisting of sequences of length $1$.
\end{proof}
\end{proposition}
Given an atomic category $\cA$ we shall call the based category obtained in this
way the \Definition{based hull} of $\cA$ and denote it $\bh{\cA}$.

Recall also (cf., \cite[p.~245]{baues01::fonct+mackey}) that for a category
$\cA$ with weak fibre products one defines a new category $\nq{\cA}$ with the
same objects as $\cA$ and where morphisms from $A$ to $B$ are formal sums of
isomorphism classes of pairs of maps $A \mapleft{f} C \mapright{g} B$ with
composition bi-additive and given by fibre products on pairs of morphisms (formal
sums are needed as the fibre product gives several maps because the fibre
products are weak). We shall refer to the basis elements as the
\Definition{primitive morphisms}. Further if $f=id$, the notation $g_*$ will be
used for the morphism and if $g=id$ the notation $f^*$ will be used. In general
the morphism equals the composite $g_*f^*$. Recall also that if $\cA$ has fibre
products then $\qh{\cA}$ denotes the category where the morphisms are
isomorphism classes of pairs $A \leftarrow C \rightarrow B$ (as sums are not
necessary). These different constructions are tightly related. For that recall
that \cite[Déf.~5.3]{baues01::fonct+mackey} defines the notion of Mackey-functor
for a category with weak fibre products (which is like the Mackey-functor of
\cite{dress73::contr} only that it is additive only over morphisms not over
objects).

However, as we are going to deal with Morita equivalences we would need to deal
with both right and left additive functors. Now,
$\nq{\cA}$ and $\dh{\cA}$ are isomorphic to their opposite categories. In fact,
the functor which is the identity on objects and takes the (class of) the
morphism $A \gets B \to C$ to $C\gets B \to A$ obtained by exchanging the
morphisms induces an equivalence between $\nq{\cA}$ and $\nq{\cA}^\circ$ resp.\
$\dh{\cA}$ and $\dh{\cA}^\circ$. This very natural anti-auto-equivalence allows
us to switch freely back and forth between right and left modules. We shall do
this without further mention. Nevertheless we should make a choice and it seems
that most of our formulas become nicer if we let a Mackey-functor be a
covariant functor with the required properties and we use the term
\Definition{left Mackey-functor} for the covariant version.
\begin{proposition}
Let $\cA$ be an atomic category.

\part $\nq{\cA}$ is equivalent to a full subcategory of $\qh{\bh{\cA}}$.

\part For any semi-additive category $\cC$ (cf.,
\cite[Déf.~4.1]{baues01::fonct+mackey}) the following data are all equivalent:
\begin{itemize}
\item A Mackey-functor on $\cA$ with values in $\cC$ in the sense of
\cite{baues01::fonct+mackey}.

\item A Mackey-functor on $\bh{\cA}$ with values in $\cC$ in the sense of
\cite{dress73::contr}.

\item An additive functor $\qh{\bh{\cA}} \to \cC$.

\item A functor $\nq{\cA} \to \cC$ that is additive on morphisms.
\end{itemize}
The correspondence is induced by the canonical functors between these categories.
\begin{proof}
The first part is clear and the second is clear using
\cite[Thm.~5.8-9]{baues01::fonct+mackey}.
\end{proof}
\end{proposition}
In the future we shall freely pass back and forth between these different ways
of looking at a Mackey-functor. In analogy with \cite[\S4]{dress73::contr} we
define for $X$ an object of $\cA$, an atomic category, $\Omega(X)$ to be the
free $\Z$-module of isomorphism classes of maps $Y \to X$. We make $\Omega$ into
a Mackey-functor by associating to \map{\alpha}{Y}{Z} the map
\map{\alpha_*}{\Omega(Y)}{\Omega(Z)} given by $[\beta] \mapsto
[\alpha\circ\beta]$ and \map{\alpha^*}{\Omega(Z)}{\Omega(Y)} by mapping
$[\beta]$ to the sum of the pullbacks of it along $\alpha$. Just as in
\cite[Prop.~4.2]{dress73::contr} any Mackey-functor $M$ has an action of $\Omega$
given by $[X \mapright{\beta} Y]a := \beta_*(\beta^*a)$ for $a \in M(Y)$.
\begin{definition}
Let $\cA$ be an atomic category. A \Definition{degree map} is a map that
associates a non-negative integer $\deg(\alpha)$ to every morphism of $\cA$
such that
\begin{itemize}
\item $\deg(\alpha\beta)=\deg(\alpha)\deg(\beta)$ for all composable morphisms
$\alpha$ and $\beta$ and

\item $\deg(\alpha) = \sum_i\deg(\alpha_i)$ where $\alpha$ is a morphism and the
$\alpha_i$ are the pullbacks of $\alpha$ along some other morphism.
\end{itemize}
\end{definition}
If $\cA$ has a degree map, $\deg$, then we may define a Mackey-functor $D$, the
\Definition{degree functor}, with $D(X)=\Z$, $\alpha_*$ is multiplication with
$\deg(\alpha)$ and $\alpha^*$ is the identity map. As always this gives an
action of $\Omega$ on $D$ which amounts to the map $\Omega \to D$ taking $X \to
Y$ to its degree which is then a map of Mackey-functors.
\begin{definition}
Let $\cA$ be an atomic category with a degree function. A
\Definition{Hecke-Mackey-functor} with values in an additive category $\cB$ is a
Mackey-functor from $\cA$ to $\cB$ for which $\beta^*\beta_*$ is multiplication
by $\deg(\beta)$ for all morphisms $\beta$ of $\cA$ or equivalently that the
action of $\Omega$ on the functor factors through the map from $\Omega$ to the
degree functor of $\cA$.
\end{definition}
\begin{remark}
\part The definition clearly does not have to make a reference to $\Omega$ and the
degree functor. The reason why we have introduced them will be clearer when we
discuss the homology Schur operad associated to a topological operad. There it
will be seen that the condition that $\beta^*\beta_*$ is multiplication by the
degree is true because we consider homology but would not be true for a
generalised homology theory. We hope to return to this point.
\end{remark}
To avoid confusion with the use of the term ``semi-additive category'' of
\cite[Déf.~4.1]{baues01::fonct+mackey} we shall call a category whose
$\Hom$-sets are given the structure of commutative monoids with composition
bi-additive \Definition{rigoids}.\footnote{So that a rigoid with one object is a
semi-ring, now often called \emph{rig}.}  One may describe also
Hecke-Mackey-functors as additive functors from a rigoid and we shall now show
how this is done. Let $A \mapleft{f} C \mapright{g} B$ be a diagram in an atomic
category. We then get a map from $C \mapright{h} D$ to one of the components of
the product of $A$ and $B$ (i.e., one of the elements of the weak fibre product
of $A$ and $B$ over the final element). We shall say that $A \leftarrow C \to B$
is a \Definition{graph} if $h$ is an isomorphism. We now define $\dh{\cA}$ to be
the category whose objects are the same as of $\cA$ and whose morphisms are
formal sums (i.e., non-negative integer linear combinations) of isomorphism
classes of graphs. The composition is the same as that in $\nq{\cA}$ with the
difference that a digram $A \mapleft{f} C \mapright{g} B$ is replaced by $A
\mapleft{p} D \mapright{q} B$ multiplied by $\deg{h}$, where $p$ and $q$ are the
projections. Again we shall call the basis elements the \Definition{primitive
morphisms}. Almost by definition we have an additive functor $\nq{\cA} \to
\dh{\cA}$ which is the identity on objects and maps $A \leftarrow C \to B$ to
$\deg{h}$ times $A \leftarrow D \to B$.
\begin{proposition}
A Mackey-functor on an atomic category with a degree function $\cA$ seen as an
additive functor on $\nq{\cA}$ is a Hecke-Mackey-functor precisely when it
factors through $\nq{\cA} \to \dh{\cA}$.
\begin{proof}
Using the notations in the preceding paragraph we have
\begin{displaymath}
g_*f^*=(qh)_*(ph)^*=q_*h_*h^*p^*
\end{displaymath}
and a Hecke-Mackey-functor maps it to
$\deg{h}\,q_*p^*$ which means that it factors through $\dh{\cA}$. On the other
hand if Mackey-functor factors through $\dh{\cA}$ then $h_*h^*$ maps to
$\deg{h}$ which means that it is a Hecke-Mackey-functor.
\end{proof}
\end{proposition}
The prototypical case of an atomic category and a degree function is obtained by
considering a finite group $G$ and a set $\cG$ of subgroups closed under
conjugation and intersection and let $\cP(G,\cG)$ be the category of transitive
permutation $G$-actions whose stabilisers are contained in $\cG$. The degree
function is given by letting $\deg(S \to T)$ be the cardinality of any fibre of
$S \to T$. Associated to this example is also a prototypical
Hecke-Mackey-functor namely the one that to a permutation set $S$ takes the
$G$-permutation representation, $\Z[S]$, on it. To a map \map{f}{S}{T} one lets
$f_*$ be the map $f_*[s] = [f(s)]$ and $f^*[t]= \sum_{s \in
f^{-1}\{t\}}[s]$. (This is also the motivation for the use of the name of Hecke
as endomorphism algebras of permutation modules will soon be seen to be
involved.) Note that as the definition of $f_*$ and $f^*$ does not use inverses
we also get a Hecke-Mackey-functor by $S \mapsto \N[S]$, where $\N[S]$ is the
free abelian monoid on $S$. 

We shall be interested essentially in one example of this construction namely
when the group is $\Sigma_n$, the group of permutations on $n$ letters, and
$\cS_n$ consists of the subgroups that fix some partition of $[1^n]$. In fact
our main interest will be in a category that is equivalent to it, namely the
category $\Mj_n$ of multisets of cardinality $n$ and multijections between
them. To any such multiset we associate the $\Sigma_n$-set $\Sigma^\alpha :=
\Hom_{\Mj}([1^n],\alpha)$ with $\Sigma_n$ acting on $[1^n]$ by the inverse of
the natural action. We denote by $\Mj$ the category of finite
multisets and multijections between them. It is the disjoint union of the
$\Mj_n$.
\begin{proposition}\label{multipermutations}
The functor $\alpha \mapsto \Sigma^\alpha$ is an equivalence of the category of
multisets of cardinality $n$ and multijections with the category
$\cP(\Sigma_n,\cS_n)$. It takes the monoidal structure of $\Mj$ to
$(S,T) \mapsto \Sigma_{m+n}\times_{\Sigma_m\times\Sigma_n}(S\times T)$ for $S \in
\cP(\Sigma_m,\cS_m)$ and $T \in
\cP(\Sigma_n,\cS_n)$. Finally, the degree function of $\cP(\Sigma_n,\cS_n)$ is
transferred to $\deg(\alpha \to \beta)={\beta\choose\alpha}$, where for a finite
multiset $\alpha$ we put
\begin{displaymath}
\alpha! := \prod_{s \in \underline{\alpha}}\mu_\alpha(s)!
\;\;\mathrm{and}\;\;{\beta\choose\alpha} =\frac{\beta!}{\alpha!}.
\end{displaymath}
\begin{proof}
We define a quasi-inverse in the following manner. Given a transitive
$\Sigma_n$-set $S$ we consider the set $\underline{T}$ of orbits of the action of $\Sigma_n$ on
$S \times [1^n]$. Each such orbit maps equivariantly onto $S$ and in particular
its cardinality is divisible by the cardinality of $S$. We then define a
multiset $T$ supported on $\underline{T}$  by $\mu_T(t) = |t|/|S|$. Given a
multiset $I$ of cardinality $n$ we get a map
$\Hom_{\Mj}([1^n],I)\times[1^n] \to \underline{I}$ that maps $(s,\varphi)$
to $\varphi(s)$ which is easily seen to induce an isomorphism of the multiset
associated to the $\Sigma_n$-set $\Hom_{\Mj}([1^n],I)$ and $I$. Conversely
given an object $S$ of $\cP(\Sigma_n,\cS_n)$ we have a map $S \to
\Hom_{\Mj}([1^n],S\times_{\Sigma_n}[1^n])$ which takes $s \in S$ to $i
\mapsto \overline{(i,s)}$ which is easily seen to be a bijection.

As for the degree function, by multiplicativity, if $\beta \to \gamma$ is a
final multijection, then $\deg(\alpha \to \beta)=\deg(\alpha \to
\gamma)/\deg(\beta \to \gamma)$ and $\deg(\alpha \to \gamma)$ by definition
equals the cardinality of $\Sigma^\alpha$ which is easily seen to equal
$|\alpha|!/\alpha!$. Finally, we deal with the monoidal structures. If
$\alpha$ and $\beta$ are multisets of cardinality $m$ and $n$ respectively then
a multijection $[1^{m+n}] \to \alpha \disjunion \beta$ gives rise to a partition
of $[1^{m+n}]$ in one part of size $m$ and another of size $n$. On the first
part we get a multijection onto $\alpha$ and on the second part one onto
$\beta$. Conversely, any such partition together with a pair of multijections
give rise to a multijection $[1^{m+n}] \to \alpha \disjunion \beta$. This makes
the description of the induced monoidal structure on $\cP(\Sigma_n,\cS_n)$
clear.
\end{proof}
\end{proposition}
Mainly with the case of $\Mj$ in view we shall mean by a \Definition{symmetric
monoidal atomic category} an atomic category with a symmetric monoidal structure
$\Tensor$ such that $A\Tensor -$ preserves weak fibre products for all objects
$A$. Furthermore, if the category also has a degree function we shall say that
it is \Definition{compatible} with the monoidal structure if $\deg(f\Tensor
g)=\deg(f)\deg(g)$.
\begin{proposition}\label{monoidal}
\part If $\cA$ is a symmetric monoidal atomic category then there is a canonical
symmetric monoidal structure on $\nq{\cA}$ making $f \mapsto f_*$ and $f \mapsto
f^*$ monoidal functors. If $\cA$ has a degree function compatible with the
monoidal structure then $\dh{\cA}$ has a symmetric monoidal structure making the
canonical functor $\nq{\cA} \to \dh{\cA}$ symmetric monoidal.

\part The disjoint union makes $\Mj$ a symmetric monoidal atomic category whose
degree function is compatible with the monoidal structure.
\begin{proof}
The first part is a simple verification left to the reader. As for
the second the preservation of weak fibre products is clear or can be inferred
from Proposition \ref{multipermutations} and the compatibility with the degree
function also follows from Proposition \ref{multipermutations}.
\end{proof}
\end{proposition}
\end{subsection}
\begin{subsection}{Almost right exact functors}

We need an extension result for functors that is well known for additive functors. For
this we need to use part of the Dold-Puppe construction (cf.,
\cite{dold61::homol+funkt}) of left-derived functors of a non-additive
functor. Recall that the left-derived functors of a functor from an abelian
category $\cA$ with sufficiently many projectives to another abelian category is
obtained as follows: For an object $A$ of $\cA$ we take a projective resolution
$P_\cdot \to A$ of it, turn it into a simplicial object $\Delta(P_\cdot)$
through the Dold-Puppe construction, apply the functor to its components,
consider the complex associated to this simplicial object and then finally take
its homology; $L_iF(A) := H_i(F(\Delta(P_\cdot)))$. 

We shall only be interested in the right-most derived functor, $L_0F(A)$ and it
only depends on the first two steps of the resolution which allows for the
following more direct formulation. An \Definition{s-display} in an abelian
category $\cA$ consists of three objects $A_1$, $A_0$, and $A$ of $\cA$, two
maps \map{\partial_0,\partial_1}{A_1}{A_0}, one map \map{s_0}{A_0}{A_1}, and one
map \map{\epsilon}{A_0}{A} with the relations
$\partial_1s_0=\partial_1s_0=\symb{id}_{P_0}$ and for which
$A_1\equaliser{\partial_1}{\partial_0}A_0\mapright{\epsilon} A$ makes $A$ the
coequaliser of $\partial_0$ and $\partial_1$. Putting $A'_1:=\ker{\partial_1}$
we get an ordinary presentation $A'_1\mapright{\partial_0}A_0\to A \to 0$ and
conversely given such a presentation we may construct an s-display by putting
$A_1 := A'_1\Dsum A_0$ and let $\partial_1$ and $s_0$ be the projection resp.\
inclusion on the last factor. This shows that the category of s-displays is
equivalent to the category of exact sequences $A'_1\mapright{\partial_0}A_0\to A
\to 0$. If $A_1$ and $A_0$ are projective we shall speak of an
\Definition{s-presentation} (of $A$).

If $F$ is a functor from $\cA$ to some other abelian category we shall say that
$F$ is \Definition{almost right exact} if for every s-display
$(A_1,A_0,A,\partial_0,\partial_1,s_0,\epsilon)$,
$F(A_1)\equaliser{F(\partial_1)}{F(\partial_0)}F(A_0)\to F(A)$ is a coequaliser
(or equivalently that
\begin{displaymath}
(F(A_1),F(A_0),F(A),F(\partial_0),F(\partial_1),F(s_0),F(\epsilon))
\end{displaymath}
is an s-display).

If $\cA$ has enough projectives and $F$ is a functor from $\cA$ to some other
abelian category we define $L_0F(A)$, for $A$ an object of $\cA$, as the
coequaliser of $F(P_1)\equaliser{\partial_1}{\partial_0}F(P_0)$ of an
s-presentation for $A$. A crucial step in the proof that the derived functors do
not depend on the resolution chosen is that homotopies between maps of complexes
correspond to homotopies between the associated simplicial objects. The notion
for s-presentations that corresponds to homotopies between maps between
presentations can not be directly formulated in terms of the s-presentations
themselves but needs the next step in the simplicial object $\Delta(P_1 \to
P_0)$. As we shall make no use of the full particulars of such a correspondence
we leave the details to the reader. For our purposes it is enough to note that a
homotopy \map{h}{P_0}{Q'_1} between maps \map{f,g}{(P'_1\to P_0)}{(Q'_1\to Q_0)} of
presentations induces a map \map{h_0}{P_0}{Q_1} for which $\partial_1h_0=f_0$
and $\partial_0h_0=g_0$ which implies that $L_0F(f)=L_0F(g)$.

We have a canonical natural transformation $L_0F \to F$ which is an isomorphism
on projective objects of $\cA$. Note furthermore that the definition of
$L_0F(A)$ gives an extension of a functor $F$ defined only on the full
subcategory of $\cA$ whose objects are the projective objects of $\cA$. The
following proposition justifies calling such an extension the \Definition{almost
right exact extension} of $F$.
\begin{proposition}
Let $F$ be a functor from an abelian category $\cA$ with sufficiently many
projectives to another abelian category. Then $F$ is almost right exact precisely when
the natural map $L_0F \to F$ is an isomorphism. In particular the almost right exact
extension of a functor defined on the full subcategory of projective objects of
$\cA$ is almost right exact.
\begin{proof}
Given an exact sequence $A'_1 \to A_0 \to A \to 0$ we may construct a
commutative diagram
\begin{displaymath}
\begin{array}{ccccccc}
P'_1&\to&P_0&\to&A&\to&0\\
\downarrow&&\downarrow&&\|\\
A'_1&\to&A_0&\to&A&\to&0
\end{array}
\end{displaymath}
with exact rows, the map $P_0 \to A_0$ surjective, and $P'_1$ and $P_0$
projective. This induces a map of displays
\begin{displaymath}
\begin{array}{ccccc}
F(P'_1\Dsum P_0)&\equaliser{}{}&F(P_0)&\to&L_0F(A)\\
\downarrow&&\downarrow&&\downarrow\\
F(A'_1\Dsum A_0)&\equaliser{}{}&F(A_0)&\to&F(A),
\end{array}
\end{displaymath}
where by assumption the upper row is a coequaliser. Assume now that $L_0F(C) \to
F(C)$ is an isomorphism for all $C$. This implies that $F(P) \to F(B)$ is
surjective for any projective object $P$ and surjective $P \to B$. Hence $F(P_0)
\to F(A_0)$ is surjective and a simple diagram chase shows that the lower row is
a coequaliser.

Conversely, if the lower row is exact for all s-displays then $F$ preserves
surjective maps so $F(P_0) \to F(A_0)$ is again surjective and then a diagram
chase concludes also in this case. The last part is then obvious.
\end{proof}
\end{proposition}
Let us now switch to the case when $\cA$ is the category of finitely generated
$R$-modules over some Noetherian ring $R$ and the target category is the
category of all $R$-modules. Starting with a functor defined on finitely
generated $R$-modules we may extend it to the category of all modules by writing
an arbitrary module $M$ as a directed limit of finitely generated modules,
$M=\dli M_\alpha$, and then put $F(M):=\dli F(M_\alpha)$. This gives a
well-defined extension as $R$ is Noetherian, the category of $R$-modules is
equivalent to the category of ind-objects of the category of finitely generated
$R$-modules. Starting with a functor defined only on the category of finitely
generated projective modules we may first consider its almost right exact
extension to finitely generated modules and then extend it to all modules. The
crucial property needed for this extension is that the category of all
$R$-modules is equivalent to the category of ind-objects of the category of
finitely generated modules (cf., \cite[App.~Prop.~2]{hartshorne66::resid}). The
two properties that is used in that proof is that every $R$-module is an
inductive limit of finitely generated $R$-modules and that
\begin{displaymath}
\dli_j\ili_i\Hom_R(G_i,H_j) = \Hom_R(\dli G_i,\dli H_j),
\end{displaymath}
where $(G_i)$ and $(H_j)$ are inductive systems of finitely generated modules.
By a slight abuse of language we shall call that final extension the
\Definition{almost right exact extension} and generally we shall say that a
functor defined on all $R$-modules is \Definition{almost right exact} if it
preserves s-displays and commutes with inductive colimits.

Finally, it is clear that a functor defined on the category of finitely
generated projective modules is determined by its restriction to the full
subcategory of finitely generated free modules and that any such functor defined
only on the subcategory has an extension to projective modules; we shall use
this fact without further mention.
\end{subsection}
\begin{subsection}{Divided power functors}

The divided power functors will play a fundamental r\^ole in what is to follow
and in this section we shall recall and establish some of its basic
properties. The quickest way of constructing them is probably as the homogeneous
components of the free divided power algebra on the module in question. Hence
recall (cf., \cite{roby63::lois}) that a divided power
algebra is a commutative augmented $R$-algebra $A$, $R$ a commutative ring, with
operators $\gamma_n$ from the augmentation ideal $I_A$ of $A$ to $A$ fulfilling
the identities of $r \mapsto r^n/n!$. The free divided power algebra $\Gamma(M)$
on an $R$-module $M$ is positively graded being $R$ in degree $0$, $M$ in degree
$1$, and the degree $n$-part $\Gamma^n(M)$, $n > 0$, generated as $R$-module by
the $\gamma_n(m)$, $m \in M = \Gamma^1(M)$. Using this we can show that
$\Gamma^n(-)$ is an almost right exact functor.
\begin{proposition}
\part The functors $\Gamma^n(-)$ are almost right exact.

\part The tensor product of almost right exact functors is almost right
exact. In particular $\Gamma^{n_1}(-)\Tensor\Gamma^{n_2}(-)\Tensor\cdots
\Tensor\Gamma^{n_k}(-)$ are almost right exact.
\begin{proof}
If $M_1 \to M_0 \to M \to 0$ is exact then $\Gamma(M_0) \to \Gamma(M)$ is
surjective and the kernel is the ideal generated by the image in $\Gamma(M_0)$
of the augmentation ideal $\Gamma^+(M_1)$ of $\Gamma(M_1)$. Indeed, the
surjectivity is clear as by universality $\Gamma(M)$ is the smallest divided
power subalgebra of itself that contains $M$ and the description of the kernel
follows from the fact that the ideal generated by the image of $\Gamma^+(M_1)$ is a
divided power ideal and hence the quotient is a divided power algebra and we
conclude by universality of $\Gamma(M)$.

Furthermore, as the tensor product of two divided power algebras has a divided
power structure compatible with the given ones we have that $\Gamma(M\Dsum
M')=\Gamma(M)\Tensor_R\Gamma(M')$. Hence, if $M'_1\Dsum M_0 \equaliser{}{} M_0
\to M$ is an s-display then we have that $\Gamma(M'_1\Dsum
M_0)=\Gamma(M'_1)\Tensor_R\Gamma(M_0)$ and hence the coequaliser of $\Gamma(M'_1\Dsum M_0)
\equaliser{}{} \Gamma(M_0)$ is $\Gamma(M_0)$ modulo $\Gamma(M_0)\Gamma^+(M'_1)$,
i.e., the ideal generated by (the image of) $\Gamma^+(M'_1)$ which, as we have
just seen is $\Gamma(M)$.

That the divided power algebra commutes with directed colimits follows directly
from the universal property.

As for tensor products of almost right exact functors being almost right exact
this is trivial for commutation with directed colimits and easy for preservation
of s-displays (and is a very trivial part of the Eilenberg-Zilber formula).
\end{proof}
\end{proposition}
\begin{remark}
\part By the same argument we get that for instance the symmetric powers are almost
right exact.

\part The almost right exactness (and in fact right exactness)
is also established in \cite[Thm.~IV.5]{roby63::lois} and is then used to show
that the kernel of $\Gamma(M_0) \to \Gamma(M)$ is the ideal generated by
$\Gamma^+(M_1)$.
\end{remark}
When the module is projective we have a well-known alternative way of describing
the divided power functors that we also shall make use of (cf.,
\cite[Prop.~IV.5]{roby63::lois}). In fact, if $P$ is a projective module, then
$\Gamma^n(P)$ can be identified with the invariants under $\Sigma_n$, the
symmetric group on $n$ letters, acting by permuting the tensor factors of
$P^{\tensor n}$. The map is given by taking $\gamma_n(p)$ to $p\tensor
p\tensor\dots\tensor p \in P^{\tensor n}$. That it is an isomorphism
$\Gamma^n(P) \to (P^{\tensor n})^{\Sigma_n}$ is most easily seen by first
reducing to $P=R^m$ for some positive integer $m$ by first writing $P$ as a
direct factor of a free module and then using commutation with direct
limits. For $R^m$ one reduces to $R$ by using behaviour of the divided power with
respect to direct sums and finally that case is trivial. With essentially the
same type of argument one also shows that
$\Gamma^{n_1}(P)\Tensor\Gamma^{n_2}(P)\Tensor\cdots \Tensor\Gamma^{n_k}(P)$ is
the module of invariants of the action of $\Sigma_{n_1}\times \Sigma_{n_2}\times
\cdots\times \Sigma_{n_k}$ on $P^{\tensor n}$, $n=\sum_in_i$.

We take the opportunity to introduce the notation 
\begin{displaymath}
\Gamma^\alpha(M) := \Gamma^{n_1}(M)\Tensor\Gamma^{n_2}(M)\Tensor\cdots
\Tensor\Gamma^{n_k}(M),
\end{displaymath}
where $\alpha$ is the \emph{ordered partition} $[n_1,n_2,\dots,n_k]$. In fact
$\Gamma^\alpha(-)$ can be canonically extended to any finite multiset. Indeed, for a
multiset $\alpha$ of cardinality $n$ we put, for a projective $R$-module $P$,
\begin{equation}\label{multi DP}
\Gamma^\alpha(P) := \left(\Dsum_{\phi \in
\Hom_{\Mj}([1^n],\alpha)}P^{\tensor n}\right)^{\Sigma_n}.
\end{equation}
We then shall also, without further mention, extend $\Gamma^\alpha(-)$ to
all $R$-modules as an almost right exact functor.

Using a slightly different way of viewing $\Gamma^\alpha(-)$ we can provide it
with a structure of Hecke-Mackey-functor.
\begin{proposition}
\part Let $G$ be a finite group and $\cG$ a set of subgroups of $G$ closed under
conjugation and intersection. Then the Hecke-Mackey-functor $S \mapsto \N[S]$ from
$\cP(G,\cG)$ to the category $S(G)$ of additive representations of $G$ on commutative
monoids induces an equivalence of categories from $\dh{\cP(G,\cG)}$ to the full
subcategory of $S(G)$ whose objects are of the form $\N[S]$ for $S \in
\P(G,\cG)$. Differently put $S \mapsto \N[S]$ is the universal
Hecke-Mackey-functor on $\cP(G,\cG)$.

\part There is a canonical structure of (left) Hecke-Mackey-functor from the atomic category
of multisets of cardinality $n$ and multijections to the category of strict
polynomial functors (on f.g.\ projective $R$-modules) which on $\alpha$ takes
the value $\Gamma^\alpha(-)$. Furthermore, that structure considered as a functor
from $\dh{\Mj}$ to the category of strict polynomial functors is monoidal, where
$\Mj$ is given the monoidal structure of proposition \ref{monoidal} and strict
polynomial functors the monoidal structure given by the tensor product.
\begin{proof}
Starting with the first part, it is enough to show that the functor
$\dh{\cP(G,\cG)}\to S(G)$ induced by $S \mapsto \N[S]$ is fully faithful. Now,
$\Hom_{G}(\N[G/H],\N[S])$ has as basis maps given by the orbits on $S$ of
$H$. On the other hand a graph $G/H \leftarrow D \to S$ gives rise to an
$H$-orbit of $S$ as the image in $S$ of the inverse image under $D \to S$ of $H$
under $G/H \leftarrow D$ and it is well-known that this map gives a bijection
between graphs and such orbits. It is equally well-known that this map maps the
basis of $\Hom_{\dh{\cP(G,\cG)}}(G/H,S)$ consisting of graphs to the basis of
$\Hom_{G}(\N[G/H],\N[S])$ given by $H$-orbits in such a way that its additive
extension realises the action of the functor $\dh{\cP(G,\cG)}\to S(G)$ on
morphisms.

As for the second part, the definition of $\Gamma^\alpha(P)$ may be rewritten as
$\Gamma^\alpha(P)=\Hom_{R[\Sigma_n]}(R[\Sigma^\alpha],P^{\tensor n})$. Hence
$\alpha \mapsto \Gamma^\alpha(-)$ factors as
\begin{displaymath}
\alpha \mapsto \N[\Sigma^\alpha] \mapsto R[\Sigma^\alpha] \mapsto
\Hom_{R[\Sigma_n]}( R[\Sigma^\alpha],(-)^{\tensor n})
\end{displaymath}
and is hence the composite of a Hecke-Mackey-functor and additive functors (one
of which is contravariant) and is thus a (left) Hecke-Mackey-functor. To show
that it is monoidal we use (\ref{multipermutations}) which gives us that
\begin{displaymath}
\Sigma^{\alpha\disjunion\beta} = \Sigma_{m+n}\times_{\Sigma_m\times\Sigma_n}(\Sigma^\alpha\times\Sigma^\beta)
\end{displaymath}
so that if $\gamma:=\alpha\disjunion\beta$ we have that $R[\Sigma^{\gamma}]$ is
the $\Sigma_{m+n}$-module induced from the $\Sigma_m\times\Sigma_n$-module
$R[\Sigma^\alpha\times\Sigma^{\beta}]$ and hence
\begin{displaymath}
\Gamma^{\gamma}(P) =
\Hom_{\Sigma_m\times\Sigma_n}(R[\Sigma^\alpha\times\Sigma^{\beta}],P^{\tensor(m+n)})=\Gamma^{\alpha}(P)\Tensor\Gamma^{\beta}(P).
\end{displaymath}
\end{proof}
\end{proposition}
We can use the two descriptions to describe certain natural transformations of
these generalised divided power functors. The most natural way of describing
them is to use the algebra structure. Hence we have a multiplication map,
\map{m_{[n_1,\dots,n_k]}}{\Gamma^{n_1}(P)\Tensor_R\Gamma^{n_1}(P)\Tensor_R\cdots\Tensor_R\Gamma^{n_k}(P)}{\Gamma^{n}(P)},
where $n := \sum_in_i$, induced by the multiplication in the divided power
algebra. On the other hand, we have a coproduct $\Gamma(P) \to
\Gamma(P)\Tensor_R\Gamma(P)$ characterised by being a map of divided power
algebras and taking $P$ (in degree $1$) to the diagonal in $P\Tensor 1\Dsum
1\Tensor P$ (also in degree $1$). It maps $\gamma_n(p)$ to
$\sum_{i+j=n}\gamma_i(p)\tensor\gamma_j(p)$. By projecting onto the
$(i,j)$-component it gives a map $\Gamma^n(P) \to
\Gamma^i(P)\Tensor\Gamma^j(P)$, for $i+j=n$, taking $\gamma_n(p)$ to
$\gamma_i(p)\tensor\gamma_j(p)$. Iterating the coproduct and projecting gives us
a map
\map{\Delta_{[n_1,\dots,n_k]}}{\Gamma^n(P)}{\Gamma^{n_1}(P)\Tensor_R\Gamma^{n_1}(P)\Tensor_R\cdots\Tensor_R\Gamma^{n_k}(P)},
when $n := \sum_in_i$.

In terms of symmetric tensors the first map, given by
$\gamma_{n_1}(p_1)\tensor\cdots\tensor\gamma_{n_k}(p_k) \mapsto
\gamma_{n_1}(p_1)\cdots\gamma_{n_k}(p_k)$, takes a tensor $m \in V^{\tensor n}$
that is invariant under
$\Sigma_{n_1}\times\Sigma_{n_2}\times\cdots\times\Sigma_{n_k}$ to the tensor
$\sum_i\sigma_im$ where the $\sigma_i$ are coset representatives for
$\Sigma_{n_1}\times\Sigma_{n_2}\times\cdots\times\Sigma_{n_k}$ in $\Sigma_n$ and
is hence invariant under $\Sigma_n$. On the other hand the map $\gamma_n(p)
\mapsto \gamma_{n_1}(p)\tensor\cdots\tensor\gamma_{n_k}(p)$ is given by the
inclusion of $\Sigma_n$-invariant tensors in the
$\Sigma_{n_1}\times\Sigma_{n_2}\times\cdots\times\Sigma_{n_k}$-invariant ones.
\begin{proposition}\label{prod/coprod}
Let $n_i$, $i=1,\dots,k$, be non-negative integers with $n:= \sum_in_i$ and let
\map{f}{[n_1,\dots,n_k]}{[n]} be the unique (final) map. Then the
(left) Hecke-Mackey-functor $\alpha \mapsto \Gamma^{\alpha}(-)$ takes $f^*$ to
$\Delta_{[n_1,\dots,n_k]}$ and $f_*$ to $m_{[n_1,\dots,n_k]}$.
\begin{proof}
To begin with, $f^*$ is induced by the $R$-linear map $R[\Sigma^\alpha] \to R$,
where $\alpha := [n_1,\dots,n_k]$, taking each basis element to $1$. This makes
it clear that the obtained map $\Gamma^n(P) \to \Gamma^\alpha(P)$ is the
inclusion of $\Sigma_n$-invariants in
$\Sigma_{n_1}\times\cdots\times\Sigma_{n_k}$-invariants. On the other hand $f_*$
is induced the $R$-linear map $R \to R[\Sigma^\alpha]$ taking $1$ to the sum of
the basis elements. If we let $e$ be a basis element fixed by
$\Sigma_{n_1}\times\Sigma_{n_2}\times\cdots\times\Sigma_{n_k}$ then the basis is
given by the elements $\sigma_ie$, where the $\sigma_i$ are coset
representatives for
$\Sigma_{n_1}\times\Sigma_{n_2}\times\cdots\times\Sigma_{n_k}$ in $\Sigma_n$,
which makes it clear that $f_*$ gives $m_{[n_1,\dots,n_k]}$.
\end{proof}
\end{proposition}
\begin{remark}
Note that every multijection is the disjoint union of final multijections and that
$\Gamma^\alpha(-)$ is monoidal as a functor in $\alpha$. As furthermore every
morphism of $\dh{\Mj}$ is of the form $f_*g^*$ this means that the proposition
can be used to give a description of the Hecke-Mackey-functor $\Gamma^\alpha(-)$.
\end{remark}
We shall need some variations on the theme of definition of $\Gamma^\alpha$. Let
$\alpha$ be a finite multiset of cardinality $n$, $P$ an f.g.\ projective
$R$-module, and \map{e}{\underline{\alpha}}{P} a function. For each multijection
\map{f}{[1^n]}{\alpha} we define an element $e_f:= e(f(1))\tensor\dots\tensor
e(f(n))$ in $P^{\tensor n}$. It is clear that $\dsum_f e_f$ gives a
$\Sigma_n$-invariant element of $\Dsum_f P^{\tensor n}$ and hence an element
$\gamma_\alpha(e) \in \Gamma^\alpha(P)$. (Note that if $\alpha =
[n_1,\dots,n_k]$ then
$\gamma_\alpha(e)=\gamma_{n_1}(e(1))\tensor\dots\tensor\gamma_{n_k}(e(k)) \in
\Gamma^{n_1}(P)\Tensor\dots\Tensor\Gamma^{n_k}(P)$.)

Further if $\alpha$ is a finite multiset of cardinality $n$ and $\{V_s\}_{s \in
\underline{\alpha}}$ is a collection of $R$-modules then we define
$\Gamma^\alpha(V)$ by almost right exactness and by the following condition if
$V$ takes projective values: For each multijection \map{f}{[1^n]}{\alpha} we
define $V^f$ by $V_{f(1)}\Tensor_RV_{f(2)}\Tensor_R\dots\Tensor_RV_{f(n)}$ and
then $\Gamma^\alpha(V)$ as the $\Sigma_n$-invariants of $\Dsum_fV^f$. (Again, if
$\alpha = [n_1,\dots,n_k]$, then $\Gamma^\alpha(V) =
\Gamma^{n_1}(V_1)\Tensor\dots\Tensor\Gamma^{n_k}(V_k)$.) Finally, if
\map{f}{\alpha}{S} is a multi-map and $\{V_s\}_{s \in S}$ is a collection of
$R$-modules, then we define $\Gamma^f(V_\bullet)$ as
$\Gamma^\gamma(V_{g(\bullet)})$, where $\gamma$ is the graph of $f$ and
\map{g}{\underline{\gamma}}{S} is the projection on the second factor. Note that
even though the second definition depends on the first we may also express the
first in terms of the second in that $\Gamma^\alpha(V)=\Gamma^f(V_\bullet)$
where $f$ is the identity function seen as a map $\alpha \to
\underline{\alpha}$. Using these notations we may express how $\Gamma^\alpha$
decomposes when applied to a direct sum.
\begin{proposition}\label{direct sums}
Let $\{V_s\}_{s \in S}$ be a collection of $R$-modules and $\alpha$ a
multiset. Then we have a natural isomorphism
\begin{displaymath}
\Gamma^{\alpha}(\Dsum_{s \in S}V_s) = \Dsum_{\map{f}{\alpha}{S}}\Gamma^f(V_\bullet),
\end{displaymath}
where $f$ runs over all multi-maps from $\alpha$ to $S$.
\begin{proof}
This is just an invariant formulation of the fact that $\Gamma^*(U\Dsum
V)=\Gamma^*(U)\Tensor\Gamma^*(V)$ and in fact the map from the right hand side
to the left hand is given by multiplication in the divided power algebra, the
factors multiplied grouped according to $f$.
\end{proof}
\end{proposition}
We have a certain functoriality for multijections: If \map{f}{\beta}{\alpha} is
a multijection and $U_t := V_{f(t)}$ then we get maps
\map{f^*}{\Gamma^\alpha(V)}{\Gamma^\beta(U)} and
\map{f_*}{\Gamma^\beta(U)}{\Gamma^\alpha(V)}, where $f^*$ and $f_*$ are defined
by the same formulas as for \map{f^*}{\Gamma^\alpha}{\Gamma^\beta} and
\map{f_*}{\Gamma^\beta}{\Gamma^\alpha}, i.e., are defined by tensor products of
the (homogeneous components of) the coproduct resp.\ the product. In particular,
if \map{g}{\alpha}{S} is a map then we get maps
\map{f^*}{\Gamma^g(V_\bullet)}{\Gamma^{g\circ f}(V_\bullet)} and
\map{f_*}{\Gamma^{g\circ f}(V_\bullet)}{\Gamma^{g}(V_\bullet)}.
\end{subsection}
\end{section}
\begin{section}{Polynomial functors}

In this section we shall introduce the notion of strict polynomial functors on
the category of $R$-modules. This notion will be almost but not quite the direct
extension of that of \cite[Def.~2.1]{fried97::cohom}. There are two reasons for this
difference. The first, and minor, is that it solves the problem that many of the
functors that will be of interest to us would not be strict polynomial functors
otherwise. This is solved in similar circumstances by introducing the notion of
analytic functors that are the direct limits of strict polynomial functors. (This
approach would in our case suffer from the slight technical problem that strict
polynomial functors are not in general determined by their underlying functors.)
The more important reason for adopting the definition that is to be presented is
that it conforms well with our eventual aim of defining the notion of polynomial monad
and in particular characterising them as essentially those monads whose algebras
admit scalar extensions (cf.\ \cite{ekedahl::schur}).

We begin by recalling the notion of polynomial maps in the sense of Roby (cf.,
\cite{roby63::lois}). If $M$ and $N$ modules over $R$, then a
\Definition{polynomial map} from $M$ to $N$ is a natural transformation
\map{f}{-\Tensor_RM}{-\Tensor_RN} of set-valued functors, where for an
$R$-module $K$, $-\Tensor_RK$ is the functor on the category of $R$-algebras
which take $S$ to $S\Tensor_RK$. As every $R$-algebra is an inductive limit of
finitely generated $R$-algebras it is easy see that we may deal only with
finitely generated $R$-algebras and we shall do that in the future. Recall
further that to $f$ we can associate unique polynomial maps $f_i$,
$i=0,1,\dots$, with the property that for each $m \in S\Tensor_RM$ and each
$\lambda \in S$, $S$ some $R$-algebra, we have that $f(\lambda m)=\sum_i
f_i(\lambda m)$ (all but a finite number of terms being $0$) and for which
$f_i(\lambda m)=\lambda^if_i(m)$. The $f_i$ are
constructed by writing $f(t\tensor m) = \sum_it^i\tensor f_i(m)$, where $t$ is
the generator of $S[t]$. We shall call $f_i$ the \Definition{$i$'th homogeneous
component} of $f$ and say that $f$ is \Definition{homogeneous of degree $i$} if
$f=f_i$. Note that (cf., \cite{roby63::lois}) homogeneous maps $M \to N$ of
degree $n$ correspond to $R$-linear maps $\Gamma^n(M) \to N$, where the identity
map $\Gamma^n(M) \to \Gamma^n(M)$ corresponds to the universal homogeneous map
$M \to \Gamma^nM$ taking $m$ to $\gamma_n(m)$.

The following rather technical lemma will be needed later.
\begin{lemma}
Let $R$ be a commutative ring and $M$ and $N$ $R$-modules. 

\part\label{decomp:i} Let \map{f}{S}{\Hom_S(S\Tensor_RM,S\Tensor_RN)} be a natural
transformation of functors in $R$-algebras $S$. There are unique $f_i \in \Hom_R(M,N)$,
$i=0,1,2,\dots$, such that for every element of $M$ all but a finite
number of the $f_i$ vanishes on it and such that for every $R$-algebra $S$ and
every $\lambda \in S$ we have that $f(\lambda) =
\sum_i\lambda^i\symb{id}_S\tensor f_i$ (the sum being finite when evaluated on
an element of $S\Tensor_RM$).

\part\label{decomp:ii} Let $Q$ be an $R$-module and $f$ a natural transformation
\map{f}{-\Tensor_RQ}{\Hom_-(-\Tensor_RM,-\Tensor_RN)} for which $f_S(\lambda q)
= \lambda^nf_S(q)$ for all $R$-algebras $S$, $\lambda \in S$, and $q \in
S\Tensor_RQ$. Then there is a unique homogeneous polynomial map
\map{F}{Q}{\Hom_R(M,N)} of degree $n$ such that for every $R$-algebra $S$ $f_S$
is the composite of \map{F_S}{S\Tensor_RQ}{S\Tensor_R\Hom_R(M,N)} and the
canonical map $S\Tensor_R\Hom_R(M,N) \to \Hom_S(S\Tensor_RM,S\Tensor_RN)$.
\begin{proof}
For the first part we consider $f(t)$, where $t$ is the variable of
$R[t]$. It is an $R[t]$-map $R[t]\Tensor_RM \to R[t]\Tensor_RN$ which is the
same thing as an $R$-map $M \to R[t]\Tensor_RN$. We can therefore define
$R$-maps \map{f_i}{M}{N} characterised by $f(t)(m)=\sum_it^i\tensor f_i(m)$
which is a finite sum as the target is $R[t]\Tensor_RN$. For any $\lambda \in S$
we have an $R$-map $R[t] \to S$ taking $t$ to $\lambda$ and so by functoriality
$f(\lambda) = \sum_i\lambda^i\symb{id}_S\tensor f_i$. 

As for the second part we choose a presentation $R[U] \to R[V] \to M \to 0$,
where $R[U]$ and $R[V]$ are the free modules on the sets $U$ and $V$. Hence for
any $R$-algebra $S$ we have an exact sequence
\begin{displaymath}
0 \to \Hom_S(S\Tensor_RM,S\Tensor_RN) \longrightarrow \prod_V S\Tensor_RN
\longrightarrow  \prod_U S\Tensor_RN
\end{displaymath}
and thus $f_S$ is given by a $V$-tuple $f^v_S$, $v \in V$, of maps $S\Tensor_RQ \to
S\Tensor_RN$ which map to a $U$-tuple with constant value $0$. Every component
$f^v_S$ is natural in $S$ and hence gives a polynomial map \map{f^v}{Q}{N} of
degree $n$. This in turn corresponds to an $R$-linear map $\Gamma^nQ \to N$ and
together they give an $R$-linear map $\Gamma^nQ \to \prod_V N$ which maps to
zero in $\prod_U N$ and thus gives a map $\Gamma^nQ \to \Hom_R(M,N)$, i.e., a
homogeneous polynomial map $Q \to \Hom_R(M,N)$ of degree $n$.
\end{proof}
\end{lemma}
We are now ready to give the definition of strict polynomial functors. Some of
our arguments will also need the stronger notion that directly correspond to
definition of \cite{fried97::cohom}.
\begin{definition}
\part A \Definition{strict polynomial functor} over a Noetherian commutative ring $R$
consists of a function $F$ from finitely generated projective $R$-modules to
$R$-modules together with the choice for each pair $P$ and $Q$ of finitely
generated projective $R$-modules and $R$-algebra $S$ a map
\begin{displaymath}
\map{F_{P,Q,S}}{\Hom_S(S\Tensor_RP,S\Tensor_RQ)}{\Hom_S(S\Tensor_RF(P),S\Tensor_RF(Q))}
\end{displaymath}
such that
\begin{enumerate}
\item If $S \to S'$ is a map of $R$-algebras and $f$ is a module map
\map{f}{S\Tensor_RP}{S\Tensor_RQ} then $F_{P,Q,S'}(id\tensor f)=id\tensor
F_{P,Q,S}(f)$ (i.e., $F$ \Definition{commutes with extension of scalars}),

\item $F_{P,P,S}(id)=id$, and

\item $F(f\circ g)=F(f)\circ F(g)$ for all $P_1$, $P_2$, $P_3$, $f \in
\Hom_S(S\Tensor_RP_2,S\Tensor_RP_3)$, and $g \in
\Hom_S(S\Tensor_RP_1,S\Tensor_RP_2)$.
\end{enumerate}

\part A \Definition{strongly polynomial functor} over a Noetherian commutative ring $R$
consists of a function $F$ from finitely generated projective $R$-modules to
$R$-modules together with the choice for each pair $P$ and $Q$ of finitely
generated projective $R$-modules and $R$-algebra $S$ a map
\begin{displaymath}
\map{F_{P,Q,S}}{\Hom_S(S\Tensor_RP,S\Tensor_RQ)} {S\Tensor_R\Hom_R(F(P),F(Q))}
\end{displaymath}
such that the relations 1)
-3) (with some obvious modifications) as for a strict
polynomial functor are fulfilled.

\part A strict (resp.\ strongly) polynomial functor $F$ is
\Definition{homogeneous of degree $n$} if
$F(\lambda\symb{id}_{S\Tensor_RP})=\lambda^n\tensor\symb{id}_{F(P)}$ for
every finitely generated projective $R$-modules, $R$-algebra $S$, and $\lambda
\in S$.
\end{definition}
Note that we have a natural map 
\begin{displaymath}
S\Tensor_R\Hom_R(F(P),F(Q)) \to \Hom_S(S\Tensor_RF(P),S\Tensor_RF(Q))
\end{displaymath}
and using it associates to any strongly polynomial functor a strict polynomial
functor.
\begin{remark}
The notion of strongly polynomial functors corresponds directly to that of
\cite[Def.~2.1]{fried97::cohom} but as the map 
$S\Tensor_R\Hom_R(F(P),F(Q)) \to \Hom_S(S\Tensor_RF(P),S\Tensor_RF(Q))$
is an isomorphism when $F(P)$ is a
finitely generated projective $R$-module we see that if $F(P)$ is finitely
generated projective when $P$ is, then a strongly polynomial structure on $F$ is
the same thing as a strict one. As this is always the case in
\cite{fried97::cohom} we see that our definition of strict polynomial functor
includes theirs.
\end{remark}
Note that if $F_\alpha$ are strongly polynomial functors then $\Dsum_\alpha
F_\alpha$ given by $(\Dsum_\alpha F_\alpha)(P):= \Dsum_\alpha F_\alpha(P)$ has
an obvious structure of strict polynomial functor which actually is a direct sum
in the categorical sense. This is not true for strongly polynomial functors
which is the main reason for introducing our definition of strict polynomial
functor.  Our first result shows that this essentially is the only difference
between strict and strongly polynomial functors.
\begin{proposition}\label{homog decomposition}
\part A strict polynomial functor $F$ can be written in a unique fashion as
$\Dsum_iF_i$ where $F_i$ is a strict polynomial functor homogeneous of degree $i$.

\part A homogeneous strict polynomial functor has a unique strongly polynomial
structure.
\begin{proof}
For a finitely generated projective $R$-module $R$ we consider the linear map $R
\to \Hom_R(P,P)$ taking $1$ to $\symb{id}_P$. Composing it with the structure
maps
\begin{displaymath}
\Hom_S(S\Tensor_RP,S\Tensor_RP) \to \Hom_S(S\Tensor_RF(P),S\Tensor_RF(P))
\end{displaymath}
gives us a situation to which (\ref{decomp:i}) applies and it gives as $u_i \in
\End_R(P)$ such that $F(\lambda)=\sum_i\lambda^iu_i$. Continuing as in
\cite[App.~A:2]{macdonald95::symmet+hall} we get that $F(st)=\sum_i(st)^iu_i$,
for $s$ and $t$ the variables of $R[s,t]$ and on the other hand $F(st)=F(s)\circ
F(t)=(\sum_is^iu_i)(\sum_it^iu_i)$ which shows that the $u_i$ are orthogonal
idempotents so that $F(P)$ splits as a direct sum of $F_i(P) := \Im(u_i)$. It is
now easily seen that for any map \map{f}{P}{Q} $F(f)$ preserves this
decomposition and hence gives a structure of strict polynomial functor on each
$F_i$ and it is equally easy to see that $F_i$ is homogeneous of degree $i$.

Now if $F$ is homogeneous of degree $n$ we have for projective $P$, $Q$, $S$ an
$R$-algebra, $f \in \Hom_S(S\Tensor_RP,S\Tensor_RQ)$ and $\lambda \in S$ that
\begin{displaymath}
F(\lambda f)=F((\lambda\symb{id}_Q)\circ f) = F(\lambda)\circ F(f) = \lambda^n F(f).
\end{displaymath}
This means that we may apply (\ref{decomp:ii}) and conclude that the structure
maps 
\begin{displaymath}
\Hom_S(S\Tensor_RP,S\Tensor_RQ) \to \Hom_S(S\Tensor_RF(P),S\Tensor_RF(Q))
\end{displaymath}
factors (uniquely) into maps 
\begin{displaymath}
\Hom_S(S\Tensor_RP,S\Tensor_RQ) \to
S\Tensor_R\Hom_R(F(P),F(Q))
\end{displaymath}
which means that we have provided $F$ with a unique
structure of strongly polynomial functor which also is homogeneous of degree
$n$.
\end{proof}
\end{proposition}
Using the proposition we may concentrate our attention on homogeneous
functors. Our current aim is to get the analogue of \cite[2.10]{fried97::cohom}
that constructs a projective generator for the category of functors of fixed
homogeneity.
\begin{proposition}\label{proj generators}
Let $R$ be a commutative Noetherian ring. 

\part The category of strict (resp.~strongly) polynomial functors over $R$ is
abelian. A sequence is exact precisely when evaluation on all finitely generated
projective modules is exact.

\part For a projective $R$-module $P$ the function $Q \mapsto
\Gamma^n\Hom_R(P,Q)$ has a structure of homogeneous strict polynomial functor of
degree $n$ and it represents the functor $F \mapsto F(P)$ on homogeneous functors
of degree $n$. In particular, $\Gamma^n\Hom_R(P,-)$ is a projective object.

\part $\Gamma^n\Hom_R(R^n,-)$ is a projective generator on the category of
homogeneous functors of degree $n$.
\begin{proof}
As both categories by Proposition \ref{homog
decomposition} are the products of the categories of homogeneous functors of fixed
homogeneity we can reduce to the homogeneous case and then again by
(\ref{homog decomposition}) to the case of homogeneous strongly polynomial
functors (of degree $n$ say).

For strongly homogeneous functors the non-trivial property is the existence of
kernel and cokernel of a map $F \to G$ together with the property that evaluated
on projective modules they are the kernel and cokernel of the evaluation of $F$
and $G$. This amounts to showing that if we define $H(P)$ as the kernel
(resp.~cokernel) of $F(P) \to G(P)$ then we want to provide it with a structure
of strongly polynomial functor (compatible with that on $F$ resp.~$G$). 

Now note that if $P$ is a finitely generated projective $R$-module then the map
$S\Tensor_R\Hom_R(P,Q) \to \Hom_S(S\Tensor_RP,S\Tensor_RQ)$ is an isomorphism so
that the structure map 
\begin{displaymath}
\Hom_S(S\Tensor_RP,S\Tensor_RQ) \to \Hom_S(S\Tensor_RF(P),S\Tensor_RF(Q))
\end{displaymath}
of a strongly polynomial functor is the same thing as a polynomial map
\begin{displaymath}
\Hom_R(P,Q) \to \Hom_R(F(P),F(Q)).
\end{displaymath}
As $F$ in our case is homogeneous of degree $n$ the map is homogeneous of degree
$n$ and hence corresponds to an $R$-linear map $\Gamma^n\Hom_R(P,Q) \to
\Hom_R(F(P),F(Q))$ (and similarly for the other functors). Hence we need to
associate to every $f \in \Gamma^n\Hom_R(P,Q)$ an $R$-linear map $H(P) \to
H(Q)$. However, by assumption we can to $f$ associate one $R$-map $F(P) \to
F(Q)$ and one $G(P) \to G(Q)$ making the diagram
\begin{displaymath}
\begin{array}{ccc}
F(P)& \to& G(P)\\
\downarrow&&\downarrow\\
F(Q)& \to& G(Q)
\end{array}
\end{displaymath}
commutative and hence inducing a unique map $H(P) \to H(Q)$ which is clearly
$R$-linear.

As for the second part, following \cite[Thm.~2.10]{fried97::cohom} we consider
the map $\Gamma^n\Hom_R(P,-) \to \Hom_R(F(P),F(-))$ which induces a map $F(P)
\to \Hom(\Gamma^n\Hom_R(P,-),F(-))$. On the other hand, we get a map
$\Hom(\Gamma^n\Hom_R(P,-),F(-)) \to F(P)$ by evaluating some
\map{f}{\Gamma^n\Hom_R(P,-)}{F(-)} on $P$ and considering the image of
$\gamma_n(\symb{id})$. These two maps are easily seen to be inverses to each
other.

As for the last part, as the category of homogeneous functors is abelian it is
enough to show that if $F$ is homogeneous of degree $n$ and $F(R^n)=0$ then
$F(P)$ vanishes for all finitely generated projective $P$. For this again it is
enough to show it when $P=R^m$ for some $m$ (as every $P$ is a factor of some
$R^m$). Furthermore we may assume that $m > n$.

If \map{e_i}{R^m}{R^m} is the projection on the $i$'th coordinate we may
consider $F(\sum_{1 \le i \le m}t_ie_i)$ for independent variables $t_i$. It is a sum
$\sum_\alpha t^\alpha u_\alpha$ and by considering the transformation $t_i
\mapsto st_i$, $s$ being a new polynomial variable, and using that $F$ is
homogeneous of degree $n$ we get that $u_\alpha=0$ unless $|\alpha|=n$. This
means that for every such $\alpha$ there is a subset $S \subset \{1,2,\dots,m\}$
of cardinality $n$ such $u_\alpha$ factors through $F$ applied to the composite $R^m \to R^n
\to R^m$ where the first map is the projection on the $S$-coordinates and the
second the inclusion. This show that all the $u_\alpha$ are zero and by putting
all the $t_i=1$ we get that $\symb{id}_{F(R^m)}=F(\symb{id})=0$.
\end{proof}
\end{proposition}
\begin{remark}
That the category of strict polynomial functors is abelian can be proved
directly from the definition when $R$ is a field as then all $R$-algebras $S$
are flat. In the general case there is a problem with kernels.
\end{remark}
In the applications we have in mind an alternative description of strict
polynomial functors will also be useful.
\begin{proposition}\label{polynomial rephrasing}
Let $R$ be a commutative Noetherian ring. A strict polynomial functor $F$
amounts to specifying for each finitely generated $R$-algebra $S$ an almost right
exact functor $F_S$ from $S$-modules to $S$-modules together with a natural
isomorphism $T\Tensor_RF_S(M) \riso F_T(T\Tensor_RM)$ for every $R$-algebra map
$S \to T$ fulfilling an evident transitivity condition. The correspondence is
then given by $F = F_R$.
\begin{proof}
Start with a strict polynomial functor $F$ and let $S$ be a finitely generated
$R$-algebra.  We first define $F_S$ on finitely generated free $S$-modules by
setting $F_S(S^n):=S\Tensor_RF(R^n)$ and then use
\begin{displaymath}
\map{F_{R^n,R^m,S}}{\Hom_S(S\Tensor_RR^n,S\Tensor_RR^m)}{\Hom_S(S\Tensor_RF(R^n),S\Tensor_RF(R^m))}
\end{displaymath}
to define $F_S$ on maps. We then extend $F_S$ uniquely to finitely generated
projective modules (as noted above) and then to all $S$-modules by almost right
exactness. By construction we have an isomorphism $T\Tensor_RF_S(M) \riso
F_T(T\Tensor_RM)$ when $M$ is finitely generated free which then extends by
additivity to an isomorphism for $M$ finitely generated projective and then to
all $M$ by almost right exactness (and right exactness of the tensor product).

The other direction is just a question of retracing the steps.
\end{proof}
\end{proposition}
This proposition allows us in particular to compose two strict polynomial
functors $F$ and $G$ by putting $(F\circ G)_S(P):=F_S(G_S(P))$.
\end{section}
\begin{section}{$\Sch$-modules}

In this section we shall make a closer study of strict polynomial functors. We
know, by the generalisation of \cite[Thm.~2.10]{fried97::cohom}, a projective
generator for functors of a fixed homogeneity which by Morita theory gives a
module description of the category. However we may (as in
\cite[Cor.~2.12]{fried97::cohom}) split up this projective generator in
components and it is more natural to consider not one but several projective
generators. Then it is also more natural to interpret Morita theory as giving an
equivalence with the category of $R$-linear functors from an $R$-linear category
with finitely many objects to the category of $R$-modules. Finally we use
previously obtained results to give a purely combinatorial description of it in
terms of multisets.

We start by introducing, for a strict polynomial functor $F$ and a multiset
$\alpha$ the $R$-module $F_\alpha := \Hom_{\cP}(\Gamma^\alpha,F)$, $\cP$ being
the category of strict polynomial functors.
\begin{proposition}
Let $F$ be a strict polynomial functor.

\part \label{decomposition} We have
\begin{displaymath}
F(\Dsum_{i\ge 1}Re_i) = \Dsum_\alpha F_\alpha,
\end{displaymath}
the sum running over all ordered partitions. Furthermore, the factor $F_\alpha$
is the image of $\gamma_\alpha(e)\tensor F_\alpha$, where $e(i)=e_i$, under the
natural map $F_\alpha\Tensor_R\Gamma^\alpha(-) \to F$.

\part The $\Gamma^\alpha$, where $\alpha$ runs over the finite multisets, form a
set of projective generators.
\begin{proof}
We start by noticing that by Proposition \ref{proj generators} and just as in the
field case (cf., \cite[Thm.\ 2.10]{fried97::cohom}) $\Gamma^n\Hom_R(R^m,-)$
represents the functor $F \mapsto F(R^m)$ on homogeneous functors of degree
$n$. Again just as in the field case (cf., \cite[Cor.\ 2.12]{fried97::cohom}) it
decomposes as the direct sum $\Dsum_\beta \Gamma^\beta(-)$, where $\beta$ runs
over all ordered partitions of $n$ supported in $[1^m]$ and that if
$\beta=[n_1,n_2,\dots,n_m]$ then $\Hom_{\sP}(\Gamma^\beta,F)$ consists of the
sub-module of $F(R^m)$ consisting of the elements of $F(R^m)$ that are
multi-homogeneous of type $[n_1,n_2,\dots,n_m]$. Furthermore, $\gamma_\alpha(e)
\in \Gamma^\alpha(\Dsum_{1\le i \le m}Re_i)$ is the universal element. More
invariantly, for a multiset $\alpha$, $\Hom_{\sP}(\Gamma^\alpha,F)$ consists of
the submodule of $F(\Dsum_{s \in \underline{\alpha}}Re_s)$ of the elements that
are $\alpha$-homogeneous and $\gamma_\beta(e) \in \Gamma^\beta(\Dsum_{s \in
\underline{\alpha}}Re_s)$ is the universal element. Now, $F(\Dsum_{i\ge 1}Re_i)$
can be decomposed into multi-homogeneous components by considering the
homogeneous decomposition of $F(\sum_it_i\pi_i)$, $\pi_i$ being the projection
onto the $Re_i$-factor which gives the decomposition $F(\Dsum_{i\ge 1}Re_i) =
\Dsum_\alpha F_\alpha$ and the fact that the $F_\alpha$ is the image of
$\gamma_\alpha(e)\tensor F_\alpha$ follows from the fact that $\gamma_\alpha(e)$
is the universal element.

To prove the second part the fact that the $\Gamma^\alpha$ are projective
follows from the first part as $F_\alpha$ is a direct factor of $F(\Dsum_{i \ge
1}Re_i)$ and hence are exact. As for faithfulness we note that every f.g.\
projective $R$-module is a factor of $\Dsum_{i \ge 1}Re_i$ so if $F(\Dsum_{i \ge
1}Re_i)=0$ then so is $F$. However as $F(\Dsum_{i \ge 1}Re_i)$ is a sum of the
$F_\alpha$ we see that if they are zero then so is $F$.
\end{proof}
\end{proposition}
The proposition implies, through standard Morita theory, that the category of
homogeneous functors of degree $n$ is equivalent with the $R$-representations of
the $R$-linear category whose objects are the multisets of cardinality $n$ and
whose $R$-module of morphisms from $\alpha$ to $\beta$ is the $R$-module
$\Hom_{\sP}(\Gamma^{\beta}(-),\Gamma^{\alpha}(-))$. This $R$-module can be given
an explicit description.
\begin{proposition}\label{HomGamma}
The $R$-module $\Hom_{\sP}(\Gamma^{\beta}(-),\Gamma^{\alpha}(-))$ has a basis
consisting of the maps $f_*g^*$, where $\beta \mapleft{g} \gamma \mapright{f}
\alpha$ is a bimultijection.
\begin{proof}
After replacing $\alpha$ and $\beta$ by isomorphic multisets we may assume that
$\alpha=[a_1,\dots,a_m]$ and $\beta=[b_1,\dots,b_n]$, where the $a_i$ and $b_j$
are positive integers summing up to the common cardinality $N$ of $\alpha$ and
$\beta$ (the proposition makes sense but is trivial when the cardinalities are
distinct). A submultiset $\gamma$ of $\alpha\times\beta$ is given by its
multiplicities $n_{ij}$ at $(i,j)$ and the condition that the projections be
multijections is equivalent to $a_i=\sum_jn_{ij}$ and $b_j=\sum_in_{ij}$ for all
$i$ and $j$.

Now, elements of $\Hom_{\sP}(\Gamma^{\beta}(-),\Gamma^{\alpha}(-))$ correspond
to elements of $\Gamma^{\alpha}(R^N))$ of multi-homo\-gene\-ous type
$[b_1,b_2,\dots,b_n]$. In $\Gamma^\ell(\sum_iRe_i)$ the product
$\gamma_{k_1}(e_1)\gamma_{k_2}(e_2)\cdots\gamma_{k_r}(e_r)$ is
multi-homo\-gene\-ous of type $[k_1,k_2,\dots,k_r]$ with the constraint that
$k_1+k_2+\cdots+k_r=\ell$ and hence the elements
$\gamma_{k_{11}}(e_1)\gamma_{k_{12}}(e_2)\cdots\gamma_{k_{1N}}(e_N)\tensor\cdots
\gamma_{k_{N1}}(e_1)\gamma_{k_{N2}}(e_2)\cdots\gamma_{k_{NN}}(e_N)$ with
$k_{i1}+k_{i2}+\cdots+k_{iN}=a_i$ and $k_{1j}+k_{2j}+\cdots+k_{Nj}=b_j$ form a
basis for $\Hom_{\sP}(\Gamma^{\beta}(-),\Gamma^{\alpha}(-))$. Hence we have a
natural bijection between the sub-multisets $\gamma$ as above and a basis for
$\Hom_{\sP}(\Gamma^{\beta}(-),\Gamma^{\alpha}(-))$.

Now, by proposition \ref{prod/coprod} and the fact that $\alpha \mapsto
\Gamma^\alpha(-)$ is monoidal it follows that under this bijection the basis
element corresponding to $\gamma$ is indeed of the form $f_*g^*$.
\end{proof}
\end{proposition}
From this we get the first main theorem. Note that as $\alpha \mapsto
\Gamma^\alpha(-)$ is a (left) Hecke-Mackey-functor we get that for a strict polynomial
functor $F$ putting $F_\alpha := \Hom(\Gamma^\alpha,F)$, the association $\alpha
\mapsto F_\alpha$ is a  Hecke-Mackey-functor.
\begin{theorem}\label{HMF2PDF}
The functor $F \mapsto \{F_\alpha\}$ gives an equivalence from the category of
strict polynomial functors (over the base ring $R$) to the category of
Hecke-Mackey-functors from $\Mj$ to the category of $R$-modules.
\begin{proof}
We know that the $\Gamma^\alpha$ form a set of projective generators so as both
categories are abelian with arbitrary direct sums it is
enough to show that the map
\begin{displaymath}
R\Tensor\Hom_{\dh{\Mj}}(\beta,\alpha) \to \Hom(\Gamma^\beta,\Gamma^\alpha)
\end{displaymath}
is a bijection. This however is exactly the content of (\ref{HomGamma}) as we
know that $\Hom_{\dh{\Mj}}(\beta,\alpha)$ has a basis consisting of the
isomorphism classes of graphs from $\beta$ to $\alpha$ and those isomorphism
classes correspond exactly to the bimultijections $\beta \to \alpha$.
\end{proof}
\end{theorem}
Because of this equivalence we shall, for $M$ a Hecke-Mackey-functor from $\Mj$
to $R$-modules and $P$ an $R$-module use the notation $M(P)$ for the value of
the corresponding functor on $P$.

Due to the importance of that special case we introduce the notation $\Sch$ for
$\dh{\Mj}$. We shall also for a commutative ring $T$ use $\Sch_T$ to denote the
scalar extension of $\Sch$ to $T$ (so that $\Hom_{\Sch_T}(\alpha,\beta)$ is a
free $T$-module on (isomorphism classes of) graphs from $\alpha$ to $\beta$).
As usual an \Definition{$\Sch$-module in an additive category $\cA$} will mean
an additive functor from $\Sch$ to $\cA$ (i.e., a Hecke-Mackey-functor). A use
of the plain ``$\Sch$-module'' will refer to an $\Sch$-module in $R$-modules.
\begin{remark}
For the purposes of this article the important category is $\Sch_{\Z}$ rather
than $\Sch$. However, we believe that it could very well be important to keep
$\Sch$ as it contains more information. More precisely, for the free commutative
monoid $\N[S]$ on a set, the set $S$ can be reconstructed from it as the
indecomposable elements. As we shall see, $\Sch$ should be compared with the rigoid
generated by $\Si$, the groupoid of finite sets and bijections, and in that case the
groupoid can be recovered from its rigoid hull.
\end{remark}
It is useful to have a more explicit way of expressing the strict polynomial
functor associated to an $\Sch$-module. The most convenient way of making such a
description is to evaluate the functor on the free $R$-module $\Dsum_{i\ge
1}Re_i$ on a countable number of generators. Any finitely generated projective
module is a direct factor of that module and hence the value of the functor on
such a module is determined by the value on $\Dsum_{i\ge 1}Re_i$ together with
knowledge of the action of endomorphisms on it.
\begin{proposition}\label{explicitation}
Let $F$ be a strict polynomial functor. 

\part  For an  $R$-module $M$ the natural map
\begin{displaymath}
\Dsum_\alpha F_\alpha\Tensor_R\Gamma^\alpha(M) \to F(M),
\end{displaymath}
where $\alpha$ runs over all ordered partitions, is surjective and its kernel is
generated as $R$-module by the relations
\begin{displaymath}
bf_*\tensor a = b\tensor f_*a,\hspace{0.5cm}a'\tensor f^*b' = a'f^*\tensor b'
\end{displaymath}
for all multijections \map{f}{\alpha}{\beta}, $a \in \Gamma^\alpha(M)$, $b \in
F_\beta$, $b' \in \Gamma^\beta(M)$, and $a' \in F_\alpha$. The kernel is also
generated by the relations
\begin{itemize}
\item For an isomorphism \map{f}{\alpha}{\beta} the relation $b\tensor f_*a =
bf_*\tensor a$ for $a \in \Gamma^\alpha(M)$ and $b \in F_\beta$.

\item For an ordered partition $[n_1,n_2,\dots,n_k]$, $n_i$ positive, the multijection
\begin{displaymath}
\map{f}{[n_1,n_2,\dots,n_k]}{[n_1+n_2,n_3,\dots,n_k]}, 
\end{displaymath}
$a \in F_{[n_1,\dots,n_k]}$, $m \in M$ and $b \in \Gamma^{[n_3,\dots,n_k]}(M)$
\begin{displaymath}
a\tensor(\Delta_{[n_1,n_2]}(\gamma_{n_1+n_2}(m))\tensor b) =
af^*\tensor(\gamma_{n_1+n_2}(m)\tensor b)
\end{displaymath}

\item For an ordered partition $[n_1,\dots,n_k]$, $n_i$ positive, the multijection
\begin{displaymath}
\map{f}{[n_1,\dots,n_k]}{[n_1+n_2,\dots,n_k]},
\end{displaymath}

$a \in F_{[n_1+n_2,\dots,n_k]}$, and $m_1,m_2 \in M$
\begin{displaymath}
a\tensor\left(m_{[n_1,n_2]}(\gamma_{n_1}(m)\tensor\gamma_{n_2}(m))\tensor
b\right)=
af_*\tensor\left(\gamma_{n_1}(m)\tensor\gamma_{n_2}(m)\tensor b\right)
\end{displaymath}
\end{itemize}
\vspace{0mm}

\part[ii] Under the identification $F(\Dsum_{n \ge 0}Re_i)=\Dsum_\alpha F_\alpha$ of
(\ref{decomposition}) we have the following identifications of $f^*$ and $f_*$
on $F_\alpha$ resp.\ $F_\beta$ for a multijection \map{f}{\alpha}{\beta} between
ordered partitions. Let, for each $i \in \underline{\beta}$, $\alpha_i$ be
$f^{-1}(i)$ with the multiplicities induced from $\alpha$.
\begin{enumerate}
\item Consider the linear map
\map{g}{\Dsum_{i \ge 1}Re_i}{\Dsum_{i \ge 1}Re_i} taking $e_j$ to $e_i$ if $j
\in \underline{\alpha_i}$ and $e_j \mapsto e_j$ if $j \notin
\underline{\alpha}$. Then we have $af^* = F(g)a$ for $a \in F_\alpha$.

\item Consider the linear map $g$, over a polynomial ring over $R$, that takes
$e_i$, $i \in \underline{\beta}$, to $\sum_{j \in \underline{\alpha_i}}s_je_j$,
where the $s_j$ are polynomial variables. Then, for $b \in F_\beta$, $bf_*$ is
the component of $F(g)b$ that is homogeneous in $s_j$ of degree
$\mu_{\alpha_i}(j)$ for all $j$.
\end{enumerate}
\begin{proof}
We begin by noting that the full subcategory $\cP$ of $\Mj$ whose objects are
the ordered partitions is equivalent to $\Mj$ and hence the category of
Hecke-Mackey-functors on $\Mj$ is equivalent to the category of
Hecke-Mackey-functors on $\cP$ so the surjectivity part of the first part
follows from the fact that the $\Gamma^\alpha$ form a set of projective
generators. We then use almost left exactness to reduce to the case when $M$ is
a f.g.\ projective $R$-module $P$. It then follows from theorem \ref{HMF2PDF}
that the kernel is generated by the relations $bf\tensor a = b\tensor fa$ for
\map{f}{\alpha}{\beta} a morphism of $\dh{\cP}$, $a \in \Gamma^\beta(P)$ and $b
\in F_\alpha$.  It is clear that is enough to let $f$ run over a set of
generators for $\dh{\cP}$ as a rigoid and the rest then follows from Proposition
\ref{generators} and (\ref{prod/coprod}).

As for the first part of \DHrefpart{ii}, let \map{A}{\Gamma^\alpha}{F} be the map
which takes $\gamma_\alpha(e)$ to $a$ (using that $\Gamma^\alpha$ represents $F
\mapsto F_\alpha$). Then $af^*$ is the image of $f^*\gamma_\beta(e)$ under $A$
and by Proposition \ref{prod/coprod} we have $f^*\gamma_\beta(e) = \tensor_{i
\in \underline{\beta}}\gamma_{\alpha_i}(e)$. This means that $f^*\gamma_\beta(e)
= \Gamma^\alpha(g)\gamma_{\alpha}(e)$ and applying $A$ we get $af^*
=A(f^*\gamma_\beta(e))=A(\Gamma^\alpha(g)\gamma_{\alpha}(e))= F(g)a$. Similarly
we have, again by (\ref{prod/coprod}),
\begin{displaymath}
f_*\gamma_\alpha(e) A= \Tensor_{i \in \underline{\beta}}\prod_{j \in \underline{\alpha_i}}\gamma_{\mu_{\alpha_i}(j)}(e_j)
\end{displaymath}
and the formula for $bf_*$ follows.
\end{proof}
\end{proposition}
\begin{subsection}{The tensor structure}

Recall that we have a monoidal structure on $\Sch$ which being additive is given
by an additive functor $\Sch\times\Sch \to \Sch$. We want to give a more precise
description of the scalar extension $\Sch \Tensor_{\Sch\times\Sch} F$ of an
$\Sch\times\Sch$-module.
\begin{proposition}\label{multi-shuffle cosets}
\part[i] Let $\gamma$, $\alpha$, and $\beta$ be finite multisets. Then every
primitive morphism in $\Sch$ from $\alpha\disjunion\beta$ to $\gamma$ factors
uniquely as $h\circ(f\disjunion g)$, where \map{f}{\alpha}{\alpha '} resp.\
\map{g}{\beta}{\beta '} are primitive morphisms in $\Sch$ and \map{h}{\alpha
'\disjunion\beta '}{\gamma} is a multi-shuffle.

\part If $M$ is an $\Sch\times\Sch$-module then for every finite multiset
$\gamma$, we have, with the map $\Sch\times\Sch
\to \Sch$ being induced by the disjoint union, the equality
\begin{displaymath}
(M\Tensor_{\Sch\times\Sch}\Sch)_\gamma = \Dsum_f f \tensor M_{(\alpha,\beta)},
\end{displaymath}
where the sum runs over all multi-shuffles \map{f}{\alpha\disjunion\beta}{\gamma}
and $f \tensor M_{(\alpha,\beta)}$ is a copy of $M_{(\alpha,\beta)}$, the notation
used only to give the action of the morphisms of $\Sch$ on it using \DHrefpart{i}.
\begin{proof}
By definition $h$ is a submultiset of $(\alpha\disjunion\beta)\times\gamma$ and
hence decomposes as a disjoint union $\delta\disjunion\epsilon$ where the
projection on the first factor maps $\delta$ resp.\ $\epsilon$ into $\alpha$
resp.\ $\beta$. Let now $\alpha '$ resp.\ $\beta '$ be the multisets supported
on the images under the second projection of $\delta$ resp.\ $\epsilon$ and with
multiplicities making the second projection induce multijection $\delta \to
\alpha '$ and $\epsilon \to \beta '$. This gives the required projection and the
uniqueness is then clear which finishes \DHrefpart{i}. 

The last part is then a direct application of the first.
\end{proof}
\end{proposition}
\begin{remark}
The result should of course be compared with the result that $(m,n)$-shuffles
give a set of coset representatives for $\Sigma_m\times\Sigma_n$ in $\Sigma_{m+n}$.
\end{remark}
Using the map $\Sch\times\Sch \to \Sch$ we get a monoidal structure on the
category of $\Sch$-modules; given two $\Sch$-modules $M$ and
$N$ we construct a functor \map{M\cdot N}{\Sch\times\Sch}{\Mod_R} by $M\cdot
N(\alpha,\beta) = M(\alpha)\Tensor_RN(\beta)$ and then put $M\boxtimes N := \Sch
\Tensor _{\Sch\times\Sch}M\cdot N$. This is easily seen to give the category of
$\Sch$-modules a symmetric monoidal structure. As we shall eventually be
interested in another monoidal structure we shall call it the \Definition{tensor
structure}. It has the following interpretation in terms of the corresponding
functors.
\begin{proposition}
Let $F$ and $G$ be strict polynomial functors and $M$ resp.\ $N$ the
corresponding $\Sch$-modules. Then $M\boxtimes N$ corresponds to the strict
polynomial functor $F\Tensor_RG$.
\begin{proof}
To define a map $M\boxtimes N \to P$, where $P$ is
the $\Sch$-module of $F\Tensor G$, it is by adjunction enough to define a map
$M(\alpha)\Tensor N(\beta) \to (F\Tensor G)_{\alpha\disjunion\beta}$ functorial
in $\alpha$ and $\beta$. This is done by noticing using that elements $f \in
M(\alpha)$ and $g \in N(\beta)$ are mappings \map{f}{\Gamma^\alpha}{F} resp.\
\map{g}{\Gamma^\beta}{G} and then map such a pair to
\begin{displaymath}
\Gamma^{\alpha\disjunion\beta} = \Gamma^\alpha\Tensor\Gamma^\beta
\mapright{f\tensor g} F\Tensor G.
\end{displaymath}
To show that the obtained map $(\Sch
\Tensor _{\Sch\times\Sch}M\cdot N)_\gamma \to (F\Tensor G)_\gamma$ is an
isomorphism we use proposition \ref{multi-shuffle cosets} where we let $\gamma$ be an
ordered partition and we want hence to show that the map
\begin{displaymath}
\Dsum_{\alpha+\beta=\gamma} F_\alpha \Tensor G_\beta \to (F \Tensor G)_\gamma
\end{displaymath}
is an isomorphism. Now, by (\ref{decomposition}), We have that $(F \Tensor
G)_\gamma$ is the $\gamma$-homogeneous part of 
\begin{displaymath}
(F\Tensor G)(\Dsum_{i \ge 1}Re_i)=F(\Dsum_{i \ge 1}Re_i)\Tensor G(\Dsum_{i \ge 1}Re_) =
(\Dsum_\alpha F_\alpha)\Tensor(\Dsum_\beta G_\beta)
\end{displaymath}
and hence
\begin{displaymath}
(F \Tensor G)_\gamma = \Dsum_{\alpha+\beta=\gamma} F_\alpha \Tensor G_\beta 
\end{displaymath}
and it is clear that the obtained morphism
\begin{displaymath}
\Dsum_{\alpha+\beta=\gamma} F_\alpha \Tensor G_\beta \to \Dsum_{\alpha+\beta=\gamma} F_\alpha \Tensor G_\beta
\end{displaymath}
is the identity map.
\end{proof}
\end{proposition}
\begin{remark}
A less computational proof runs as follows: We have a pair of adjoint functors
between strict polynomial bifunctors and strict polynomial functors $F \mapsto
(U \mapsto F(U,U))$ resp.\ $G \mapsto ((U,V) \mapsto G(U\Dsum V))$. On the other
hand $\Sch\times\Sch$-modules classify strict polynomial bifunctors and it can
be seen, as has been noted in the remark after proposition \ref{explicitation},
that $G \mapsto ((U,V) \mapsto G(U\Dsum V))$ correspond to restriction of
modules along $\Sch\times\Sch \to \Sch$ and hence by adjunction $F \mapsto (U
\mapsto F(U,U))$ corresponds to scalar extension. Finally, for strict polynomial
functors $F$ and $G$, the bifunctor $(U,V) \mapsto F(U)\Tensor G(V)$ corresponds
to the $\Sch\times\Sch$-module $M\cdot N$, where $M$ and $N$ are the
$\Sch$-modules corresponding to $F$ resp.\ $G$. From this the result follows.
\end{remark}
\end{subsection}
\end{section}
\begin{section}{The $p$-local case}

Macdonald's theorem (cf. \cite[App.\ A:5.3]{macdonald95::symmet+hall}) gives
over a field of characteristic zero another description of strict polynomial
functors. We shall now show that his result follows from ours. The proof will
also be seen to have implications when rather than containing $\Q$ the ring $R$
is only $p$-local. To prepare for it we let,
for a prime $p$, $\dMj{p}$ be the full subcategory of $\dh{\Mj}$ whose objects
are those multisets all of whose multiplicities are powers of $p$ with scalars
extended to the local ring $\Z_{(p)}$. Similarly we
let  $\dMj{0}$ be the full subcategory of $\dh{\Mj}$ whose objects
are those multisets all of whose multiplicities are equal to $1$ with scalars
extended to $\Q$.
\begin{proposition}\label{p-reduction}
\part Let $p$ be a prime and $\cA$ an additive category in which multiplication
by any integer prime to $p$ on $\Hom$-groups is bijective. Then restriction and
extension of functors induce an equivalence from the category of
Hecke-Mackey-functors with values in $\cA$ to the category of additive functors
from $\dMj p$ to $\cA$.

\part Let $\cA$ be an additive category in which multiplication by any non-zero
integer on $\Hom$-groups is bijective. Then restriction and extension of
functors induce an equivalence from the category of Hecke-Mackey-functors with
values in $\cA$ to the category of additive functors from $\dMj 0$ to $\cA$.
\begin{proof}
For the first part it is enough to show that each finite multiset is a direct
factor in $\dMj p$ of a multiset all of whose multiplicities are powers of
$p$. This is done by induction over the number of multiplicities that
aren't. Hence assume that $m$ is the multiplicity at $s \in \underline{\alpha}$
and that it is not a power of $p$. Note that if \map{f}{\alpha '}{\alpha} is a
multijection that is an isomorphism outside of $s$ and if $m_1,\dots,m_k$ are
the multiplicities of the points above $s$, then $\deg f$ equals the multinomial
coefficient $m\choose m_1,\dots,m_k$ and hence $\alpha$ is a direct factor of
$\alpha '$ if that coefficient is prime to $p$. Now write $m$ in base $p$;
$a_0+a_1p+\cdots+a_np^n$ and choose $f$ such that
\begin{displaymath}
m_1,\dots,m_k =
\overbrace{\vphantom{p^n}1,1,\dots,1}^{a_0\;\mathrm{times}},
\overbrace{\vphantom{p^n}p,p,\dots,p}^{a_1\;\mathrm{times}},\dots,
\overbrace{p^n,p^n,\dots,p^n}^{a_n\;\mathrm{times}}.
\end{displaymath}
It is then well-known that the multinomial coefficient is prime to $p$.

As for the second part, we do a similar thing, only now we just need for the
multinomial coefficient to be non-zero and we can hence choose all the $m_i$ to
be equal to $1$.
\end{proof}
\end{proposition}
\begin{remark}
Note that $\dMj p$ uses all multisets for the definition of $\Hom$-sets. This is
not really necessary: Suppose $\alpha \leftarrow \gamma \to \beta$ is a
bimultijection where $\alpha,\beta \in \dMj p$ but $\gamma$ isn't. Then we may
find as in the proof of the proposition a $\gamma ' \to \gamma$ such that
$\gamma ' \in \dMj p$ and the degree $d$ of $\gamma ' \to \gamma$ is not
divisible by $p$. Thus we have $\alpha \leftarrow \gamma \to \beta = 1/d(\alpha
\leftarrow \gamma ' \to \beta) \in \Hom_{\dMj p}(\alpha,\beta)$.
\end{remark}
\begin{example}
In the $p$-local case the first case that differs from Macdonald's theorem is
the case of functors of degree $p$. The proposition then gives that the category
of such functors is equivalent to the category of tuples $(M,N,f,g)$, where $M$
is an $R[\Sigma_p]$-module, $N$ an $R$-module and \map{f}{M}{N} and
\map{g}{N}{M} are equivariant $R$-linear maps (for the trivial action of
$\Sigma_p$ on $N$) fulfilling $fg=p$ and $gf=\sum_{\sigma \in \Sigma_p}\sigma$.
\end{example}
\end{section}
\begin{section}{Scalar extension}

One way of getting a strict polynomial functor is to start with a collection
$\{M_n\}$ of right $R[\Sigma_n]$-modules and then put $S(M)(P)=\sum_nM_n\Tensor_{\Sigma_n}P^{\tensor
n}$. In that case $S(M)_\alpha =
R[\Sigma^\alpha]\Tensor_{\Sigma_n}M_n$, where $n := |\alpha|$. In particular
when $\alpha = [n_1,\dots,n_k]$ then $S(M)_\alpha =
(M_n)_{\Sigma_{n_1}\times\cdots\times\Sigma_{n_k}}$. In this section we shall
give a direct construction of the $\Sch$-module associated to $S(M)$.

If \map{h}{\cA}{\cB} is an additive  functor of rigoids, $\cC$ an
abelian category with arbitrary colimits and \map{F}{\cA}{\cC} an additive
functor we shall denote by $\cB \Tensor_\cA F$ the left Kan extension of $F$
along $h$. Note that the value of this extension on an object $b \in
\cB$ is the colimit of $F(a)$ over the category of maps $h(a) \to b$.

To begin with we shall apply this construction to the following situation: We
let $\Si$ be the groupoid of finite sets and bijections and let $\Sigma$ be its
rigoid hull, i.e., $\Hom_\Sigma(S,T) = \N[\Hom_{\Si}(S,T)]$. We have a functor
$\Si \to \Mj$ taking a set to the corresponding multiset and it extends to an
additive functor $\Sigma \to \Sch$.
\begin{proposition}
Assume that \map{M}{\Sigma}{\Mod_R} is an additive right functor. Then for any
finite multiset $\alpha$ of cardinality $n$ the value $(M\Tensor_\Sigma\Sch
)_\alpha$ of $M\Tensor_\Sigma \Sch$ at $\alpha$ equals
$M_{[1^n]}\Tensor_{R[\Sigma_n]}R[\Sigma^\alpha]$. Furthermore, the strict
polynomial functor $M(-)\Tensor_\Sigma\Sch$ associated to $M \Tensor_\Sigma
\Sch$ equals $S(M)$.
\begin{proof}
The first part follows almost directly from the description of the Kan extension
recalled above: We have that the set of maps $[1^n] \to \alpha$ is given by
$\N[\Sigma^\alpha]$ and it is easily seen that the colimit is exactly
$M_{[1^n]}\Tensor_{R[\Sigma_n]}R[\Sigma^\alpha]$.

For the rest of the first part, let $N$ be the $\Sch$-module corresponding to $
M_n\Tensor_{\Sigma_n}\Dsum_{n \ge 0}P^{\tensor n}$.  To define a map $M
\Tensor_\Sigma \Sch \to N$ we need by adjunction only define a map $M \to N$ of
$\Sigma$-modules. However, $N_n$ consists of the $[1^n]$-homogeneous elements of
$M_n\Tensor_{\Sigma_n}\Dsum_{n \ge 0}(\Dsum_{1\le i \le n}Re_i)^{\tensor n}$ and
we map $m \in M_n$ to the residue of $m\tensor e_1\tensor\cdots\tensor e_n$
which is clearly an isomorphism $M_n \to N_n$. Considering a general finite
multiset $\alpha$ we have that $N_\alpha$ is the $\alpha$-homogeneous part of
$M_n\Tensor_{\Sigma_n}\Dsum_{n \ge 0}(\Dsum_{i \in
\underline{\alpha}}Re_i)^{\tensor n}$. Such a part exists only in the summand
with $n = |\alpha|$ and then the $\alpha$-homogeneous basis vectors of
$(\Dsum_{i \in \underline{\alpha}}Re_i)^{\tensor n}$ is in $\Sigma_n$-invariant
bijection with $\Sigma^\alpha$ and we get that $N_\alpha =
M_n\Tensor_{\Sigma_n}R[\Sigma^\alpha]$. As the quotient maps $M_n \to M_\alpha$
and $N_n \to N_\alpha$ are given by the action of an element of $\Sch$, they and
the maps $M_n \to N_n$ and $M_\alpha \to N_\alpha$ form a commutative diagram
showing that $M_\alpha \to N_\alpha$ is an isomorphism.
\end{proof}
\end{proposition}
It follows from the proposition that the scalar extension from $\Sigma$ to
$\Sch$ is not exact. That means that given a $\Sigma$-module $M$ we get not only
one $\Sch$-module by scalar extension but a sequence of them;
$\Tor^{\Sigma}_*(\Sch,M) := H_*(M\Tensor_{\Sigma}^L\Sch)$. We then get the
following result as corollary.
\begin{corollary}
If $M$ is a $\Sigma$-module then for each finite multiset $\alpha$ we have
$\Tor^{\Sigma}_*(M,\Sch)_\alpha = \Tor_{\Sigma_n}(M_{[1^n]},R[\Sigma^\alpha])$.
\begin{proof}
This follows immediately from the proposition applied to a $\Sigma$-projective
resolution of $M$.
\end{proof}
\end{corollary}
A particularly interesting case is when $M = R[\Sigma_n/\rho(G)]$ for a faithful
permutation representation \map{\rho}{G}{\Sigma_n} in which case
$H_*(\Sigma_\alpha,M)$ is the sum over the $G$-orbits of $\Sigma^\alpha$ of the
homology (with coefficients in $R$) of the $G$-stabiliser of a point of the
orbit. We shall denote this homology $H^\rho_*(G)$. In \cite{ekedahl::schur} we
shall make a closer study of it and in particular show how it can be used to
express the homology of wreath products.
\end{section}
\bibliography{preamble,abbrevs,alggeom,algtop,algebra,ekedahl}
\bibliographystyle{pretex}

\end{document}